\documentclass{amsart}

\usepackage{amssymb}
\usepackage{amscd}
\usepackage{amsmath}
\usepackage{amsfonts}
\usepackage{a4}
\usepackage{amsthm}
\usepackage{color}

\theoremstyle{plain}

\newtheorem{theo}{Theorem}[section] 
\newtheorem{prop}[theo]{Proposition}
\newtheorem{lemme}[theo]{Lemme}
\newtheorem{cor}[theo]{Corollary}
\newtheorem{defin}[theo]{Definition}

\theoremstyle{definition}

\newtheorem{rem}[theo]{Remarque}

\newcommand{\Ci}{C^{\infty}}

\newcommand{\R}{{\mathbb{R}}}
\newcommand{\N}{{\mathbb{N}}}
\newcommand{\C}{{\mathbb{C}}}
\newcommand{\Z}{{\mathbb{Z}}}

\newcommand{\de}{{\delta}}
\newcommand{\be}{{\beta}}
\newcommand{\al}{{\alpha}}

\newcommand{\De}{{\Delta}}
\newcommand{\si}{{\sigma}}
\newcommand{\hb}{{\hbar}}

\newcommand{\ga}{{\gamma}}
\newcommand{\te}{{\theta}}
\newcommand{\om}{{\omega}}
\newcommand{\Om}{\Omega}
\newcommand{\ph}{{\varphi}}

\newcommand{\Ga}{{\Gamma}}

\newcommand{\La}{{\Lambda}}

\newcommand{\identite}{{\operatorname{Id}}}

\newcommand{\trace}{{\operatorname{tr}}}

\newcommand{\Hilbert}{{\mathcal{H}}}

\newcommand{\Toeplitz}{{\mathcal{T}}}
\newcommand{\Lie}{{\mathcal{L}}}

\newcommand{\K}{{{\mathcal{K}}}}
\newcommand{\Fourier}{{{\mathcal{F}}}}
\newcommand{\Comp}{\mathcal{J}}
\newcommand{\Compint}{\Comp_{\operatorname{int}}}
\newcommand{\demint}{{\mathcal{D}}_{\operatorname{int}}}
\newcommand{\demi}{\mathcal{D}}

\newcommand{\Quant}{{\mathcal{Q}}^{\operatorname{m}}}

\newcommand{\Um}{U^{\operatorname{m}}}
\newcommand{\Vm}{V^{\operatorname{m}}}
\newcommand{\compm}{\circ_{\operatorname{m}}}
\newcommand{\Groupm}{{\mathcal{G}}_0^{\operatorname{m}}}

\newcommand{\Mor}{\mathcal{M}}
\newcommand{\End}{{\operatorname{End}}}
\newcommand{\Hom}{{\operatorname{Hom}}}

\newcommand{\Preq}{{\mathcal{P}}^{\operatorname{m}}}

\newcommand{\Sym}{\operatorname{Sym}}

\newcommand{\Group}{{\mathcal{G}}}

\newcommand{\bi}{\zeta^{\frac{1}{2}}}
\newcommand{\bic}{\zeta}

\newcommand{\Op}{\operatorname{Op}}

\newcommand{\algebre}{{\mathcal{A}}}

\newcommand{\qs}{\mathcal{Q}}
\newcommand{\pqs}{\mathcal{P}}

\newcommand{\scratch}[1]{}

\author{L. Charles}

\address{Universit{\'e} Pierre et Marie Curie-Paris6, UMR 7586 Institut de
  Math{\'e}matiques de Jussieu, Paris, F-75005 France.}

\email{charles@math.jussieu.fr}

\title[Semi-classical properties of geometric quantization]{Semi-classical properties of geometric quantization with
  metaplectic correction}

\keywords{Geometric Quantization, Toeplitz operator, Fourier integral
  operator, Half-form bundle}

\subjclass{53D22, 53D50, 53D55, 81S30, 47L80}

\begin{document}

\begin{abstract} 
The geometric quantization of a symplectic manifold
endowed with a prequantum bundle and a metaplectic structure is defined
by means of an integrable complex structure. We prove that its
semi-classical limit does not depend on the choice of the complex
structure. We show this in two ways. First, by introducing unitary
identifications between the quantum spaces associated to the various
complex polarizations and second, by defining an asymptotically flat
connection in the bundle of quantum spaces over the space of complex
structures. Furthermore Berezin-Toeplitz operators are intertwined by
these identifications and have principal and subprincipal
symbols defined independently of the complex structure. The relation
with Schr{\"o}dinger equation and the group of prequantum bundle
automorphisms is considered as well. 
\end{abstract}

\bibliographystyle{plain}
\maketitle


\section{introduction} 

Geometric Quantization of Kostant \cite{Ko} and Souriau \cite{So} is a procedure which
associates a quantum space to a symplectic manifold endowed with a
prequantum bundle and a polarization. Since its introduction, there has been some attempt to find natural
{\em identifications} between quantum spaces associated to different
polarizations (cf. \cite{Bl}, \cite{Ra}). In the case of symplectic
compact manifolds with complex polarizations, Ginzburg
and Montgomery observed in \cite{GiMo} that a natural identification
does not exist for a broad class of manifolds. Recently Foth and Uribe
\cite{FoUr} obtained semi-classical results in the same direction. 

We prove that there exists a natural {\bf semi-classical}
identification when the definition of the quantum spaces is altered with the
metaplectic correction. This result is a consequence of our study
undertaken in \cite{firstpart} of the symbolic calculus of Toeplitz
operators and Lagrangian sections that we extend in this paper to Fourier integral operators.  
Before we state our results, let us discuss quantization without metaplectic correction.

\subsection{Ordinary quantization} \label{sec:intw}

Let $(M, \om)$ be a symplectic compact manifold with a prequantization
bundle $L\rightarrow M$, i.e. a Hermitian line bundle
with a connection of curvature $\frac{1}{i} \om$. Denote by
$\Compint$ the space of integrable complex structure of $M$
compatible with $\om$ and positive. To any $j \in \Compint$ is associated a sequence of quantum
spaces 
$$ \qs_k (j) := \{ j\text{-holomorphic sections of } L^k \}, \qquad k
=1,2,...$$ 
Here the holomorphic structure of $L$ is the one compatible with the
connection. The semi-classical limit corresponds to $k \rightarrow \infty$. 
When $k$ is suf\-fi\-cien\-tly large, the Kodaira vanishing
theorem and Riemann-Roch-Hirzebruch theorem imply that the dimension
of $\qs_k(j)$ is given by a Riemann-Roch number, which only depends on
the symplectic structure of $M$ and $k$. Assume that we can choose such an
integer $k$ independently of the complex structure $j$. Then for any $j_a, j_b \in \Compint$ we can
identify $\qs_k (j_a) $ with $\qs_k (j_b)$ by means of a unitary map
$$U_k(j_a,j_b) : \qs_k (j_a) \rightarrow \qs_k (j_b).$$ These
identifications are mutually compatible if they satisfy:
\begin{itemize} 
\item (functoriality) $U_k
 (j_b, j_c) \circ U_k (j_a, j_b ) = U_k (j_a, j_c )$, for any $j_a, j_b, j_c \in \Compint$.
\end{itemize}
Moreover if these maps are canonical in the sense that they only
depend on the complex, symplectic and prequantum structures, they should satisfy:
\begin{itemize} 
\item (naturality) for any prequantization bundle automorphism $\Phi$
  of $L^k$ and complex structures $j_a, j_b \in \Compint$, the diagram
$$ \begin{CD} 
 \qs_k (j_a)  @> U_k (j_a, j_b) >> \qs_k (j_b) \\
 @V{\Phi^*}VV @VV{\Phi^*}V \\
\qs_k (\Phi^* j_a)  @> U_k (\Phi^*j_a, \Phi^*j_b) >> \qs_k (\Phi^* j_b)
 \end{CD} $$
commutes. 
\end{itemize} 
Here the vertical maps are pull-back by $\Phi$, sending a
$j$-holomorphic section into a section holomorphic with
respect to $\Phi^* j:= \phi^* j$, where $\phi$ is the
symplectomorphism of $M$ covered by $\Phi$.  Sometimes one only requires
an identification between the projectivised quantum
spaces.

It is important to observe that if there exists such a collection
$\{ U_k(j_a,j_b)$, $(j_a, j_b) \in \Compint^2\}$ which is
both functorial and natural, then for any complex structure $j \in
\Compint$, the quantum space $\qs_k (j)$ becomes a representation of
the group $\Group$ of prequantization bundle automorphism of
$L^k$. Indeed let us set 
\begin{gather} \label{eq:rep} 
V_k(\Phi) :=  U_k( \Phi^* j, j) \circ \Phi^* : \qs_k (j)  \rightarrow \qs_k(j), \qquad \Phi \in \Group. \end{gather}
Then for any prequantization bundle automorphisms $\Phi_1$ and $\Phi_2$,
we have
\begin{xalignat*}{2}
V_k(\Phi_1) \circ V_k (\Phi_2)= & U_k(\Phi_1^* j,j ) \circ
\Phi_1^* \circ U_k(\Phi_2^* j ,j) \circ \Phi_2^* \\
= & U_k(  \Phi_1^* j ,j)\circ U_k( \Phi_1^* \Phi_2^* j ,\Phi_1^* j) \circ \Phi_1^* \circ \Phi_2^* \quad \text{ by naturality,} \\
= & U_k(  \Phi_1^* \Phi_2^* j ,j) \circ
 (\Phi_2\circ \Phi_1)^* \quad \text{ by functoriality,} \\
 = & V_k ( \Phi_2\circ \Phi_1 ).
\end{xalignat*} 
Considering the associated infinitesimal representation, Ginzburg and
Montgomery proved in \cite{GiMo} that the existence of such a
representation contradicts "no go" theorems in many
cases. Indeed one can view the Lie
algebra of $\Group$ as $\Ci (M, \R)$, the Lie bracket being the
Poisson bracket. Then assuming that the maps $U_k(j_a,j_b)$ depend
smoothly on $j_a$ and $j_b$, we obtain a Lie
algebra representation 
\begin{gather*} 
  \Ci (M, \R)  \rightarrow  \End ( \qs_k (j) ).\end{gather*} 
By \cite{GiMo}, since $M$ is compact and $\qs_k(j)$ is finite
dimensional, the associated
projective representation is trivial. But for a broad class of
manifolds $M$, $\Group$ contains a finite dimensional subgroup which preserves
a complex structure $j$ and whose induced representation on $ \qs_k
(j)$ is not projectively
trivial. The same arguments contradict also the existence of an
identification between the projectivised quantum spaces.

In spite of this result, there is a natural identification of a
particular interest which has been introduced for the quantization of
the moduli spaces of flat connections (cf. \cite{Hi} and \cite{AxWi}). To define it consider the quantum spaces $\qs _k (j)$
as the fibers of a bundle $\qs _k \rightarrow \Compint.$ Then introducing a functorial and
natural family $(U_k(j_a, j_b))$ which depends smoo\-th\-ly on $j_a$ and
$j_b$ amounts to endowing this bundle with a flat
$\Group$-invariant connection. Now consider $\qs_k$ as a subbundle of
$$\pqs_k := \Ci (M, L^k) \times \Compint \rightarrow \Compint$$ 
Since $\pqs_k$ is trivial, it
has a natural flat connection and $\qs_k$ is equipped with the
projected connection. 
Because of the previous result, the
curvature $R_k$ of $\qs_k$ can not vanish in general. On the
other hand by the theory of Boutet de Monvel and
Guillemin \cite{BoGu}, the Toeplitz operators provide an asymptotic  representation of the Poisson algebra
$\Ci (M)$ as operators on $\End (\qs_k (j))$ when $k \rightarrow
\infty$. So it is possible that the
curvature $R_k$ is asymptotically flat (cf. end of section \ref{sec:fdr} for a quantitative argument). Foth and Uribe compute the asymptotics of
$R_k$ in \cite{FoUr} and prove the following: for any $j \in \Compint$ and tangent vectors $\eta, \mu
  \in T_j \Compint$, there exists a function $ f(\eta, \mu) \in \Ci (M)$
such that
$$ R_k (\eta, \mu) = \Pi_k(j) f(\eta, \mu)  + O(k^{-1}) : \qs_k (j)
\rightarrow \qs_k(j)$$
where $\Pi_k(j)$ is the orthogonal projector of $\Ci(M, L^k)$ onto
$\Quant_k(j)$.
\scratch{and the $O(k^{-1})$ is for the uniform norm (pas dans
la version de \cite{FoUr}).}
Furthermore, they give a simple formula for the multiplicator $f( \eta, \mu)$, which
shows that it does not vanish for a generic choice of $(\eta,
\mu)$. Consequently the curvature is not asymptotically flat. Neither
is it asymptotically projectively flat.

\subsection{Main results} 

Let us turn to geometric quantization with metaplectic
correction. The metaplectic structures were introduced by Kostant in
\cite{Ko2} as metaplectic principal bundles lifting the symplectic frame
bundle (cf. also \cite{GuSt} and \cite{BaWe}). Here we use the half-form 
bundle approach (cf. \cite{Wo}) more convenient for our purpose. 

Given a complex structure $j \in \Compint$, a half-form bundle
$(\delta, \varphi)$ of $(M,j)$ is a line bundle $\delta \rightarrow M$ together
with an isomorphism of line bundles $$\varphi : \delta^2 \rightarrow \Lambda^{n,0}_j
T^*M $$ covering the identity of $M$. $(M,j)$ admits a half-form
bundle if and only if the second Stiefel-Whitney class of $M$
vanishes. 
From now on, we assume this condition is satisfied and we set 
$$ \Quant_k (j, \delta, \varphi) := \{ j\text{-holomorphic sections of
  $\delta\otimes L^k$} \}, \quad k=1,2,... $$ 
where the superscript "m" stands for metaplectic. Here the holomorphic structure of $\delta$ is such that $\varphi$ is an
  isomorphism of holomorphic bundles.  

There is an obvious notion of isomorphism between two
  half-form bundles associated to the same complex structure and these
  isomorphisms give rise to isomorphisms between the associated quantum
  spaces. Our aim is to extend this to the whole collection $\demi$ of triples $(j, \delta, \varphi)$, where $j$ ranges through
  $\Compint$.

In section \ref{sec:hfb}, we define a
  collection $\Mor$ of morphisms, which makes $(\demi, \Mor)$ a
  category such that every morphism is an isomorphism. Important facts
  are that the
  automorphism group of any $a \in \demi$ is $\Z_2$ and the isomorphism
  classes are in one to one correspondence with the elements of
  $H^{1}(M, \Z_2)$. Furthermore isomorphism classes correspond to
  equivalence classes of metaplectic structures of $M$.

\begin{theo} \label{theo:1} 
There exists a family $( (\Um_k (\Psi) )_k ; \;\Psi \in \Mor  )$ such
that for any morphism $\Psi: a \rightarrow b$, the sequence $(\Um_k
(\Psi))_k$ consists of operators 
$$ \Um_k (\Psi) : \Quant_k (a) \rightarrow \Quant_k (b)$$
that are unitary if $k$ is sufficiently large. Furthermore, for any composable morphisms $\Psi$ and $\Psi' \in \Mor$, we have 
$$ \Um_k (\Psi') \circ \Um_k (\Psi) = \Um_k (\Psi' \circ \Psi) + O(k^{-1})$$ 
where the estimate $O(k^{-1})$ is for the uniform norm of operators. 
\end{theo}

One of the original motivations to introduce the metaplectic correction was to define some natural pairings between the quantum
spaces associated to different polarizations, which are called now Blattner-Kostant-Sternberg pairings. Our construction of the
operators $U_k(\Psi)$ is rather different. These are Fourier integral
operators with a prescribed principal symbol and the functoriality
property is a consequence of the symbolic calculus. 

We interpret this theorem as a semi-classical functoriality of
quantization with half-form bundle. Moreover, the family $((\Um_k(
\Psi))_k;\; \Psi \in \Mor)$ is natural with respect to a
suitable action of the group $\Group$ of prequantization bundle 
automorphisms of $L$ on $(\demi, \Mor)$. We can therefore adapt the
previous construction \eqref{eq:rep} and we obtain for any $a \in \demi$ an
asymptotic representation on $\Quant_k (a)$ of a central extension by
$\Z_2$ of the identity component of $\Group$.  This
is in some sense a generalization of the standard metaplectic representation. 

We will also prove that the operators $\Um_k (\Psi)$ can be defined as
parallel transport in an appropriate bundle. Let us consider a smooth
family  $((\delta_j, \varphi_j), \; j\in \Compint)$ of isomorphic
half-form bundles. Let $\Quant_k
\rightarrow \Compint$ be the quantum space bundle, whose fiber over $j$ is the space
of $j$-holomorphic sections of $L^k \otimes \delta_j$. 

\begin{theo} \label{theo:2} 
For any positive integer $k$, the bundle $(\Quant_k \rightarrow
\Compint)$ has a canonical connection $\nabla^{\Quant_k}$.  The
sequence $( \nabla^{\Quant_k}, \; k\in \N^*)$ satisfies
\begin{itemize} 
\item for any $j \in \Compint$ and tangent vectors $\eta, \mu \in T_j
  \Compint$, the uniform norm of the curvature $R^{\Quant_k} (\eta, \mu)$ is
  $O(k^{-1})$,
\item the parallel transport in $\Quant_k$
  along a curve $\gamma$ with endpoints $j_a$ and $j_b$ is equal to
  $\Um_k (\Psi)$ modulo $O(k^{-1})$. Here $\Um_k (\Psi)$ is the sequence
  of theorem \ref{theo:1} and $\Psi$ is
  the 
  half-form bundle morphism $$(j_a, \delta_{j_a}, \varphi_{j_a})
  \rightarrow  (j_b, \delta_{j_b}, \varphi_{j_b})$$ obtained  by
  extending continuously the identity of $(j_a, \delta_{j_a},
  \varphi_{j_a})$ in morphisms $(j_a, \delta_{j_a}, \varphi_{j_a})
  \rightarrow  (\ga(t), \delta_{\ga(t)}, \varphi_{\ga(t)})$. 
\end{itemize} 
\end{theo} 
The connection $\nabla^{\Quant_k}$ is induced by a connection on the prequantum
space bundle. However the latter bundle is not trivial contrary to the case
without metaplectic correction.

The paper is organized as follows. Section 2 is devoted to
preliminary material. Section \ref{FIO} contains our results about
symbolic calculus for Fourier integral operators. These results
are reformulated in section \ref{reformulation} with the half-form bundle
formalism. In section \ref{sec:hfq}, we deduce theorem
\ref{theo:1} and related facts on the representation of the
prequantization bundle automorphisms and the Schr{\"o}dinger equation. The
study of the quantum space bundle and its connection is in sections
\ref{sec:geo}. Section \ref{sec:acG} is devoted to the action of the
prequantization bundle automorphism group on the quantum space bundle.

\section{Preliminaries} \label{section:prel} 

Let $(E, \om)$ be a symplectic real vector space of dimension
$2n$. Let $\Comp (E, \om)$
be the space of complex structures $j$ of $E$ compatible with $\om$ and
positive. Given $j \in \Comp(M, \om)$, we denote by $\La^{n,0}_{j}
E^*$ the line of complex linear forms of $E$ of type $(n,0)$ for the
complex structure $j$. 
\begin{defin} \label{def:Psi}
Given $j_a$ and $j_b \in \Comp (E,\om)$, let 
$\Psi_{j_a,j_b}$ be the linear map from $\La^{n,0}_{j_a} E^*$
to $\La^{n,0}_{j_b} E^*$ such that 
$$ \Psi_{j_a,j_b} (\al ) \wedge \bar{\beta}   = \al \wedge
\bar{\be}, \quad \forall \ \al, \beta \in \La^{n,0}_{j_a} E^*.  $$
\end{defin}

Let us give some elementary properties of these maps. First, $\Psi_{j_a,j_b}$ is well-defined and invertible because the sesquilinear pairing $$ \La^{n,0}_{j_b} E^* \times  \La^{n,0}_{j_a} E^* \rightarrow \C, \quad
(\al_b , \al_a) \rightarrow  \al_b \wedge \bar{\al}_a / \om^{\wedge n} $$
is non-degenerate, $j_a$ and $j_b$ being positive. Whenever  $j_a
=j_b$, $\Psi_{j_a,j_b}$ is the identity. With the usual scalar
product on $\La^{n,0}_{j} E^*$ defined by means of $\om$ and $j$, the
adjoint of $\Psi_{j_a,j_b}$ is $\Psi_{j_b,j_a}$. This is easily checked
using that the scalar product of $\al, \beta \in \La^{n,0}_{j} E^*$ is
given by 
$$ i ^{n(2-n)} \al \wedge \bar{\be} / \om^n .$$
Last definition that we need is the following. 

\begin{defin} Given $j_a, j_b, j_c$ in $\Comp (E, \om)$, let $\bic (j_a,j_b,j_c)$ be the
  complex number such that 
\begin{gather}   \notag
  \Psi_{
j_a,j_c} = \bic (j_a,j_b,j_c) \; \Psi_{j_b,j_c} \circ \Psi_{j_a, j_b},
\end{gather}
\end{defin} 

As we will see in the
next section, the symbols of Fourier
integral operators behave in part as square roots of the $\Psi_{j_a,j_b}$. This will appear first via the continuous
square root $\bi$  of the complex function $\bic : \Comp^3 (E, \om) \rightarrow \C^*$ determined by $\bi (j,j,j) = 1$, for any
$j \in \Comp(M, \om)$. $\bi$ is well-defined and analytic because $ \Comp^3 (E,
\om)$ is contractible and $\bic$ is an analytic function (cf. \eqref{eq:zeta}).

It follows from the associativity of the composition that $\bic$ is a
cocycle
\begin{gather} \label{eq:cocycle}
 \bic (j_b,j_c,j_d) .\bic^{-1} (j_a, j_c ,j_d ).\bic (j_a, j_b, j_d)
. \bic^{-1} (j_a,j_b,j_c ) =1 \end{gather} 
Furthermore $\Psi_{j,j}= \identite$ implies
\begin{gather} \label{eq:ident} 
 \bic (j_a,j_b,j_b ) =\bic (j_a,j_a,j_b ) =1 .\end{gather} 
The function $\bi$ satisfies the same equations. 

\scratch{Ne pas croire que $\bic ( j_a,j_b,j_a )$ est aussi {\'e}gal {\`a} $1$. On
sait seulement qu'il est r{\'e}el. Plus g{\'e}n{\'e}ralement $\Psi_{j_a,j_b}
^* = \Psi_{j_b,j_a}$ implique que 
$ \bar{\bic} ( j_a,j_b,j_c ) = \bic ( j_c,j_b,j_a)$. }

To prepare further proofs, we compute the function $\bic$
in the following  parame\-trization of $\Comp (E, \om)$. Let us choose
a fixed complex structure $j_0\in \Comp (E, \om)$. Then for any $j \in \Comp (E, \om)$,
the space of linear forms of type $(1,0)$ with respect to $j$, viewed
as a subspace of $$E^* \otimes \C = \La^{1,0}_{j_0} E^* \oplus
\La^{0,1}_{j_0} E^*,$$ is the graph of a complex linear map
\begin{gather} \label{eq:defmu}
 \mu : \La^{1,0}_{j_0} E^* \rightarrow \La^{0,1}_{j_0} E^* .\end{gather}
The condition that $j$ is compatible with $\om$ is that 
$$ \om ( \mu^t X, Y) + \om ( X, \mu^t Y) = 0, \qquad \forall \ X, Y
\in E^{0,1}_{j_0} $$
where $\mu^t : E^{0,1}_{j_0} \rightarrow E^{1,0}_{j_0}$ is the
transposed of $\mu$. 
And the positivity of $j$ translates into the positivity of the
Hermitian map:
$$ \identite - \mu^t \bar{\mu}^t :E^{1,0}_{j_0} \rightarrow E^{1,0}_{j_0}$$ 
This defines a one-to-one correspondence between $\Comp(M, \om)$ and an open set of a subspace of
$\Hom (\La^{1,0}_{j_0} E^*, \La^{0,1}_{j_0} E^*)$. 

\scratch{Introduisons une base $\te^i$ de $\La^{1,0}_{j_0} E^*$ telle que $\om = i \sum \te^i
\wedge \bar{\te}^i $.
Alors $\mu$ v{\'e}rifie les conditions pr{\'e}c{\'e}dentes si et seulement
si sa matrice $U$ donn{\'e}e par 
 $$\mu ( \te^i) = \sum_{j} U^i_j \bar{\te}^j,$$
  est sym{\'e}trique et telle que $1 - U \bar{U}$ soit d{\'e}finie positive. En dimension 1,
  on retrouve donc le disque unit{\'e}.}

For any $j$, let us identify the $(n,0)$-forms with respect to $j$ with the
$(n,0)$-forms with respect to $j_0$ by
the map 
$$  \La^{n,0}_{j_0} E^* \rightarrow \La^{n,0}_{j} E^*, \qquad
\al^1 \wedge ...\wedge \al^n \rightarrow (\al^1 + \mu (
\al^1))\wedge...\wedge (\al^n + \mu (\al^n)). $$  
Then straightforward computations prove the following lemma.

\begin{lemme} \label{lem:the_Lemme}
With the previous identifications, $\Psi_{j_a,j_b}$ regarded as a map from $
\La^{n,0}_{j_0} E^*$ to itself is the multiplication by 
$$  \det \left( \begin{array}{cc} \identite & \bar{\mu}_a \\
    \mu_a & \identite \end{array} \right) 
. \operatorname{det}^{-1}  \left( \begin{array}{cc} \identite & \bar{\mu}_a \\
    \mu_b & \identite \end{array} \right) $$
where the matrices represent maps from $\La^{1,0}_{j_0} E^* \oplus
    \La^{0,1}_{j_0} E^*$  to itself. Consequently,
\begin{gather} \label{eq:zeta}
 \bic ( j_a, j_b, j_c) = \frac{ \det \left( \begin{array}{cc} \identite & \bar{\mu}_a \\
    \mu_b & \identite \end{array} \right) 
. \operatorname{det}  \left( \begin{array}{cc} \identite & \bar{\mu}_b \\
    \mu_c & \identite \end{array} \right) }
{\det \left( \begin{array}{cc} \identite & \bar{\mu}_a \\
    \mu_c & \identite \end{array} \right) 
. \operatorname{det}  \left( \begin{array}{cc} \identite & \bar{\mu}_b \\
    \mu_b & \identite \end{array} \right)
} \ . \end{gather} 
\end{lemme} 

\scratch{En fait l'application qui me sert {\`a} identifier
  $\La^{n,0}_{j_0} E^*$ et $\La^{n,0}_{j} E^*$ n'est rien d'autre que
  $\Psi_{j_0,j}$. Donc {\`a} strictement parler, le lemme donne  
$$\Psi_{j_0,j_b}^{-1}
  \circ \Psi_{j_a, j_b} \circ \Psi_{j_0, j_a}.$$
}

\section{Fourier integral operator} \label{FIO}

Let $(M,\om)$ be a symplectic compact connected manifold with a prequantization
bundle $(L, \nabla)$, i.e. $L$ is a Hermitian line bundle and
$\nabla$ a
connection of curvature $\frac{1}{i} \om$. The quantizations of $(M,\om)$ we
will consider depend on two additional datas: a complex structure $j$ of $M$
compatible with $\om$ and positive, and a holomorphic Hermitian line
bundle $K$ over the complex manifold $(M,j)$.  

Let us denote by $\K$ the collection of such pairs $(j,K)$. To any $a =
(j_a, K_a) \in \K$, we
associate the sequence of Hilbert spaces 
$$\Hilbert_k (a) := \{  \text{ holomorphic sections of } L^k \otimes
K_a \}, \qquad k =1,2,...$$
where the holomorphic structure of $L$ is the one compatible with the
connection $\nabla$ such that $L \rightarrow M$  is holomorphic with
respect to $j_a$. The scalar product is defined by means of the Hermitian structure of $ L^k \otimes K_a$
and the Liouville measure of $M$. 

For any $a,b \in \K$, let us introduce the space $\Fourier ( a ,
b)$ of Fourier integral operators  from $\Hilbert (a)$ to $\Hilbert ( b)$. Their definition is a slight generalization of the one
in \cite{oim2} because of the fiber bundles $K_a$ and $K_b$. Consider a sequence
$(S_k)$ such that for every $k$, $S_k$ is an operator $\Hilbert _k
(a) \rightarrow \Hilbert _k (b)$.  The scalar product of $\Hilbert _k
(a)$ gives us an isomorphism 
$$ \Hom (\Hilbert _k (a) , \Hilbert _k (b) ) \simeq  \Hilbert _k (b)
\otimes \overline{\Hilbert}
_k (a) .$$
The latter space can be regarded as the space of holomorphic sections of  
$$ (L^k \otimes K_b ) \boxtimes (\bar{L}^k \otimes \overline{K}_a)
\rightarrow M^2,$$
where $M^2$ is endowed with the complex structure $(j_b, -j_a)$. The
section $S_k (x,y)$ associated in this way to $S_k$ is its Schwartz kernel. 

We
say that $(S_k)$
is a Fourier integral operator of $\Fourier ( a ,
b)$ if 
\begin{gather} \label{def:FIO}
 S_k(x,y) =  \Bigl( \frac{k}{2\pi} \Bigr)^{n} E^k(x,y) f(x,y,k) + O
(k^{-\infty}) \end{gather}
where 
\begin{itemize}
\item 
$E$ is a section of  $L \boxtimes \bar{L}
\rightarrow M^2$ such that  $\| E(x,y)
\| <1 $ if $x \neq y$, 
$$ E (x,x) = u \otimes \bar{u}, \quad \forall u \in L_x \text{ such that }
  \| u \| = 1, $$
and $ \bar{\partial} E \equiv 0 $
modulo a section vanishing to any order along the diagonal.
\item
  $f(.,k)$ is a sequence of sections of  $ K_b  \boxtimes  \bar{K}_a
\rightarrow M^2$ which admits an asymptotic expansion in the $\Ci$
  topology of the form 
$$ f(.,k) = f_0 + k^{-1} f_1 + k^{-2} f_2 + ...$$
whose coefficients satisfy $\bar{\partial} f_i \equiv 0  $
modulo a section vanishing to any order along the diagonal.
\end{itemize}

Let us define the principal symbol of $(S_k)$ to be the map $x \rightarrow
f_0(x,x)$. Using the Hermitian structure of $K_a$, we regard it as a
section of $\Hom (K_a , K_b) \rightarrow M$. The principal symbol map 
$$\sigma :
\Fourier ( a ,b) \rightarrow \Ci(M, \Hom (K_a , K_b))$$
satisfies the expected property.

\begin{theo} \label{P0}
The following sequence is exact
$$ 0 \rightarrow \Fourier (a,b) \cap O(k^{-1}) \rightarrow 
\Fourier ( a ,b) \xrightarrow{\sigma} \Ci(M,\Hom (K_a , K_b)) \rightarrow 0, $$
where the $O(k^{-1})$ is for the uniform norm of operators. 
\end{theo}

The composition of these operators is also as expected, with some
complications regarding the product of the symbols. Given three
complex structures $j_a$, $j_b$ and $j_c$ of $M$, we denote by $\bi
(j_a, j_b, j_c )$ the function of $\Ci (M)$ whose values at $x$ is the
complex number $\bi (j_{a}(x), j_{b}(x), j_{c}(x) )$ defined in section
\ref{section:prel} with $E = T_x M$.

\begin{theo} \label{P1}
Let $a$, $b$ and $c$ belong to $\K$. If $T \in \Fourier (a,b)$ and $S
\in \Fourier (b,c)$, then $S\circ T$ is a Fourier integral operator of
$\Fourier (a,c)$. Furthermore,  
$$ \si ( S \circ T ) = \bi (j_a, j_b, j_c ) \; \si (S) \circ \si (T).
$$
\end{theo}

The two previous theorems were essentially proved in the chapter 4.1 of
\cite{oim2} except the formula for the composition of the symbols, which
will be proved in chapter \ref{sec:preuve_P1}. Since the composition
of operators is associative, the same holds for the symbol. Observe
that this can
be directly checked with the cocycle relation \eqref{eq:cocycle}.

$\Fourier (a,a)$ is  the space $\Toeplitz
(a)$ of Toeplitz operators of 
$\Hilbert (a)$. Equivalently, a Toeplitz operators is any  sequence $(T_k: \Hilbert_k(a) \rightarrow
\Hilbert_k(a) )$ of operators of the form 
$$ T_k = \Pi_k g(.,k) + R_k , $$
where $\Pi_k$ is the orthogonal projector of $L^2 (M, L^k \otimes
K_a)$ onto $\Hilbert_k (a)$, $g(.,k)$ is a sequence of $\Ci( M)$ with an asymptotic expansion
$g_0 + k^{-1} g_1 +...$ in the $\Ci$ topology and $R_k$ is  $O(k^{-\infty})$. As a result the principal symbol
$\si (T_k)$ is the function $g_0$. Let us define the {\em normalized symbol} of
$(T_k)$ to be the formal series 
$$  g(.,\hb) + \frac{\hb}{2} \Delta g(.,\hb), $$ 
where $g(.,\hb) = g_0 + \hb g_1
+ ...$
and $\Delta$ is the holomorphic Laplacian. We are actually only interested in
the two first terms of this series, which are the principal symbol and
the {\em subprincipal symbol} 
$g_1 + \tfrac{1}{2} \Delta g_0$. As a consequence of the works of Boutet
de Monvel and Guillemin, the product of the normalized symbol is a
star-product (\cite{BoGu}). 

\begin{theo}  \label{P2}
Let $S$ be a Fourier integral operator of $\Fourier (a,b)$ and 
  $T_a \in \Toeplitz (a)$,  $T_b \in
  \Toeplitz (b)$ be two Toeplitz operators with the same principal
  symbol $f$. Then 
\begin{gather*} 
 T_b \circ S - S \circ T_a = \tfrac{1}{k} R \end{gather*}
with  $R\in \Fourier (a,b)$. Furthermore the principal symbol of $R$ is 
$$ \si (R) = \bigl( f_{1,b} - f_{1,a} + \tfrac{1}{2} \langle \alpha_{j_a,j_b}
, X_f \rangle \bigr) \si (S)   + \tfrac{1}{i}
\nabla_{X_f}^{\Hom(K_a,K_b)} \si (S)
$$
where 
\begin{itemize}
\item $f_{1,a}$, $f_{1,b}$ are the subprincipal symbols of $T_a$
  and $T_b$ respectively, 
\item $X_f$ is the Hamiltonian vector field of $f$,
\item $\alpha_{j_a,j_b}$ is the one-form of $M$ such that
  $$\nabla \Psi_{j_a,j_b} = \tfrac{1}{i} \alpha_{j_a,j_b} \otimes
  \Psi_{j_a,j_b},$$
where $\Psi_{j_a,j_b}$ is the section of $$\Hom (
  \La^{n,0}_{j_a}T^*M,\La^{n,0}_{j_b} T^*M ) \rightarrow M$$ whose
  value at $x$ is the endomorphism $\Psi_{j_a(x),j_b(x)}$ defined in
  \ref{def:Psi} with $E = T_x M$. The connection $\nabla$ is induced by the Chern
  connections of $ \La^{n,0}_{j_a}T^*M$ and $\La^{n,0}_{j_b} T^*M$.
\end{itemize}
\end{theo}

By theorem \ref{P1}, $T_b S - S T_a$ is Fourier integral operator of
$\Fourier (a,b)$. Since $$\bi (j_a, j_a, j_b ) = \bi (j_a, j_b, j_b)
=1,$$ its principal symbol vanishes and consequently $R\in \Fourier (a,b)$. So the proof of the
theorem consists in computing the principal symbol of $R$. This is
postponed to chapter \ref{preuve1}. Let us deduce some interesting consequences.
Applying the theorem with two
Toeplitz operators $S$ and $T$ of $\Toeplitz
(a)$, we recover that the principal symbol of $\tfrac{k}{i}  [ S, T]$ is the
Poisson bracket of the principal symbols of $S$ and $T$. Actually, we
can also compute the subprincipal symbol of $\tfrac{k}{i}  [ S, T]$ with the previous theorem.

First, the operators of $\Fourier (a,b)$ may be used to identify  $\Hilbert
(a) $ with $\Hilbert (b)$ in a semi-classical sense. More precisely, we
consider the  space $U \Fourier (a,b)$ consisting of Fourier integral
operators of $\Fourier (a,b)$ which are unitary in the sense that
$$ S_kS^*_k =
 \identite_{\Hilbert_k (a)} \text{ and } 
 S_k^*S_k =
\identite _{\Hilbert_k (b)} $$ 
when $k$ is sufficiently large. By some standard
 argument that we briefly recall now, $ U \Fourier (a,b)$ is not empty if and
 only if $K_a$ and $K_b$ are isomorphic as line bundles. First
 it follows directly from the definition of a Fourier integral
 operator that the adjoint of an operator $S \in \Fourier (a,b)$ belongs to
 $\Fourier (b,a)$ and its principal symbol is the adjoint of
the principal symbol of $S$. So if $S \in U \Fourier (a,b)$, theorem \ref{P1} implies that
 the principal symbol of $S$ is a line bundle
 isomorphism $K_a \rightarrow K_b$.  Conversely
 if  $K_a$ and $K_b$ are
 isomorphic, there exists an elliptic  $R
 \in \Fourier (a,b)$, meaning that its principal symbol doesn't
 vanish anywhere. Then $R ^* R$ is an elliptic Toeplitz operator by theorem
  \ref{P1}. So $(R ^* R) ^{-\frac{1}{2}} $ is a Toeplitz operator
 (cf. as instance \cite{oim1}). Finally $R(R ^* R) ^{-\frac{1}{2}}$
 belongs to $U \Fourier (a,b)$.

Now if $S \in U \Fourier (a,b)$ and $T_a$ is a Toeplitz operator of
$\Hilbert (a)$, then by theorem \ref{P1} $$ T_b = S T_a
S^*$$ is a Toeplitz operator of $\Hilbert (b)$ with the same principal
symbol as $T_a$. Applying theorem \ref{P2}, we compute its
subprincipal symbol in terms of the principal and subprincipal symbols
of $T_a$:
\begin{gather} \label{trans_symb} 
 f_{1,b} = f_{1,a} + \langle \alpha_{S} -  \tfrac{1}{2} \alpha_{j_a,j_b}
, X_{f} \rangle  \end{gather}
where $  \alpha_S$ is such that  
$$
  \nabla^{\Hom(K_a,K_b)} \si(S) = \tfrac{1}{i} \alpha_{S}
\otimes \si(S).$$
A consequence of this formula is the following result. 
\begin{theo} \label{SP_norm}
The composition law $*_a$ of the normalized symbols of the  
  Toeplitz operators of $\Hilbert (a)$ satisfies:
$$ f *_a g = fg + \tfrac{\hb}{2i} \{ f, g \} + O(\hb^2)$$ 
and 
$$ \tfrac{i}{\hb} ( f *_a g - g *_a f) =  \{ f, g \} - \hb \langle
\om_{j_a} - \tfrac{1}{2} \om_{K_a} , X_f \wedge X_g \rangle  + O(\hb^2) $$
where $\frac{1}{i} \om_{j_a}$ and $\frac{1}{i} \om_{K_a}$ are the
Chern curvatures of  $\La^{n,0}_{j_a} T^*M$  and $K_a$ respectively.
\end{theo} 
Indeed, if $K_a$ and $K_b$ are isomorphic and $*_a$ satisfies the result, the same holds for $*_b$
because of \eqref{trans_symb} and the relations
$$ \om_{j_b} = d \al_{j_a, j_b} + \om_{j_a}, \quad \om_{K_b} = d
\al_{S} + \om_{K_a} .$$
Furthermore we can explicitly compute $*_a$, in the case where $M$ is $\C^n$ with
$\Hilbert (a)$ the Bargmann space, and the result is satisfied. Of
course, this is not sufficient to conclude. But it appears in the
proofs of the previous theorems that
all the results about the symbolic calculus are completely local and we can
really deduce in this way theorem \ref{SP_norm}. 

\section{Proof of theorem \ref{P1}} \label{sec:preuve_P1}

The proof relies on the complex stationary phase lemma. We only
sketch the first part, because the details appeared in \cite{oim2}, with
some typos however. 
The Schwartz kernel of an operator $S \in \Fourier (a,b)$
is by definition of the form
$$ \Bigl( \frac{k}{2 \pi} \Bigr)^n  E^{k}_{a,b} (x,y) f(x,y,k) +
O(k^{-\infty})$$
Let us write on a neighborhood of the diagonal 
$$\nabla^{L \boxtimes \bar{L}} E_{a,b} = \tfrac{1}{i}
\al_{a,b} \otimes E_{a,b} $$
The following lemma is proved in \cite{oim2}. 
\begin{lemme} The one-form $\al_{a,b}$ vanishes along the diagonal of
  $M^2$. Furthermore, for every vector fields $X_1,X_2,Y_1,Y_2$ of $M$
$$ \Lie_{(X_1,X_2)} \langle \al_{a,b}, (Y_1,Y_2) \rangle (x,x)= \om (
\bar{q}_{a,b} (X_1 - X_2) , Y_1)(x) + \om (q_{b,a} (X_1 - X_2) , Y_2)(x)$$
where $\bar{q}_{a,b}(x)$ and $q_{b,a}(x)$ are
respectively the projections onto $(T_x M)_{j_b}^{0,1}$ with kernel
$(T_xM)_{j_a}^{1,0}$ and onto $(T_x M)_{j_a}^{1,0}$ with kernel $(T_x M)_{j_b}^{0,1}$.
\end{lemme}
Consider now $S \in  \Fourier (a,b)$ and $S' \in \Fourier (b,c)$. The
Schwartz kernel of $S'S$ is
$$ \Bigl( \frac{k}{2 \pi } \Bigr)^{2n} \int_{M} E^k_{b,c}(x,y). E^k_{a,b}(y,z) f'(x,y,k) .f(y,z,k)
\ \mu_M(y) + O(k^{-\infty})$$
with $\mu_M$ the Liouville form of $M$. 
Since $|E_{b,c}(x,y).E_{a,b}(y,z)| < 1$ outside the diagonal of  $M^{3}$, this
integral is $O(k^{-\infty})$ outside the diagonal of $M^2$, and to
estimate it on a neighborhood of $(x,z) = (x,x)$ it suffices to integrate on a neighborhood of $x$. We evaluate the
result by applying the stationary phase lemma. Let us write
\begin{gather} \label{eq:def_Phi}
 E_{b,c}(x,y). E_{a,b}(y,z) = e^{i \Phi (x,y,z)} t(x)
\otimes \bar{t}(z) \end{gather}
with $t$ a unitary local section of $L \rightarrow M$. We deduce from
the previous lemma the following facts.
\begin{itemize} 
\item  $ d_{y} \Phi$ vanishes along the diagonal $\De_3$ of $M^3$. 
\item if $Y_1$ and $Y_2$ are two
tangent vectors of $M$ at $x$,
\begin{gather} \label{eq:der_sec_Phi}
 d^2_y \Phi (Y_1,Y_2) (x,x,x) = \om ( q_{c,b} Y_1, Y_2 )(x) - \om (\bar{q}_{a,b}Y_1,
Y_2)(x),\end{gather} 
In particular $d^2_y \Phi$ is non-degenerate along $\De_3$.
\item the kernel of the tangent map to 
$d_y \Phi$ at $(x,x,x)$ is 
$$ \bigl( T_{(x,x,x)} \De_3 \otimes \C \bigr) \oplus \bigl(
(T_x M)^{0,1}_{j_c}  \times (0) \times (0) \bigr) \oplus \bigl(
 (0) \times (0) \times (T_x M)^{1,0}_{j_a}  \bigr). $$
\end{itemize} 

These ensure that we can apply the stationary phase lemma (cf. \cite{Ho} or the appendix of \cite{firstpart}). Thus the
Schwartz kernel of 
$S'S$ is of the form
$$ \Bigl( \frac{k}{2 \pi} \Bigr)^n  F ^k(x,z) g(x,z,k) +
O(k^{-\infty})$$
where $g(.,k)$ is a sequence of sections of $K_c\boxtimes \bar{K}_a$
which admits an asymptotic expansion in negative power of $k$ and   
$$ F(x,z) = e^{i \Phi^r (x,z)} t(x) \otimes \bar{t}(z)$$
with
\begin{gather} \label{eq:phase_red} 
 \Phi^r (x,z) \equiv \Phi (x,y,z) \end{gather}
modulo a linear combination with $\Ci$ coefficient of the functions
$\partial_{y^i} \Phi (x,y,z), i=1,...,2n$. 

Let us check that the
section $F$ satisfies the assumptions following equation \eqref{def:FIO}. Since
$\partial_{y^i}\Phi$ vanishes along the diagonal, it follows from
\eqref{eq:def_Phi} that $\Phi^r(x,x) =0$. Furthermore, we have 

\begin{lemme} \label{lem:comp_phas}
Consider $M^2$ as a complex manifold with complex
  structure $(j_c, -j_a)$. Then $\bar{\partial} F \equiv 0$
  modulo a section vanishing to any order along the diagonal. 
\end{lemme}

\begin{proof} 
Introduce complex coordinates  $x^1,...,x^n$ on $M$ for 
$j_c$. Let us write 
 $$\nabla t = \tfrac{1}{i} \ t \otimes  \sum a_j dx^j - \bar{a}_j d\bar{x}^j .$$
Derivating equation \eqref{eq:def_Phi} and using that $
\nabla_{(\partial_{\bar{x}^i},0)} E_{b,c} $ vanishes to any order
along the diagonal of $M^2$, we get
\begin{gather} \label{eq:der_Phi}
 \partial_{\bar{x}^i} \Phi (x,y,z) + \bar{a}_i(x) \equiv 0
\end{gather} 
modulo $\mathcal{I}_{\Delta_3} (\infty)$, i.e. modulo a function
vanishing to any order along the diagonal of $M^3$. 
Thus $$\partial_{\bar{x}^i} \partial_{y^j} \Phi
(x,y,z) \equiv 0 \mod \mathcal{I}_{\Delta_3} (\infty).$$ In the
same way, if $z^1,...,z^n$ are complex coordinates for  $j_a$, we show
that   
$$\partial_{z^i} \partial_{y^j} \Phi
(x,y,z) \equiv 0 \mod \mathcal{I}_{\Delta_3}
(\infty).$$ Then we deduce from \eqref{eq:phase_red} and \eqref{eq:der_Phi} that for any multi-index $\al$
and $\be$, the function 
$$  \partial_{\bar{x}^1} ^{\al(1)} ...\partial_{\bar{x}^n}^{\al(n)} \partial_{z^1} ^{\be(1)} ...\partial_{z^n} ^{\be(n)}
\bigl( \partial_{\bar{x}^i} \Phi^r(x,z) + \bar{a}_i (x) \bigr)$$
vanishes along the diagonal $\De_2$ of $M^2$. 
This implies that 
$$ \nabla^{L \boxtimes \bar{L}}_{(\partial_{\bar{x}^i},0)} F \equiv 0$$
modulo a section vanishing to any order along $\De_2$. We treat
in the same way the covariant derivatives of $F$ with respect to 
the vector fields $(0,\partial_{z^i})$.
\end{proof} 

Then since the kernel of $S'S$ is a holomorphic section of  $$(L^k \otimes K_b) \boxtimes (\bar{L}^k \boxtimes K_a)$$
the coefficients of the asymptotic expansion of $g(.,k)$
satisfy $\bar{\partial} g_l \equiv 0$ modulo a section vanishing at
any order along the diagonal. So we proved that $S'S$ is a Fourier
integral operator of $\Fourier (a,c)$. 

Final step is to compute its
symbol. By the stationary phase lemma, we have 
$$ g(x,x,k) = f'(x,x,k).f''(x,x,k) \frac{ \delta(x)}{ \operatorname{det}^{\frac{1}{2}} [ -i \partial_{y^j} \partial_{y^k} \Phi
  (x,x,x)]_{j,k}} + O(k^{-1}) $$
where  $$\mu_M (y) = \delta (y) .|dy^1...dy^{2n}|$$
We deduce from  \eqref{eq:der_sec_Phi} that
\begin{xalignat*}{2}  
 -i d_y^2 \Phi (Y_1, Y_2) (x,x,x) & = - \om ( i q_{c,b} Y_1 - i \bar{q}_{a,b}Y_1,
Y_2) (x)\\ &=  - \om ( j_b q_{c,b} Y_1 + j_b \bar{q}_{a,b}Y_1,
Y_2) (x)\\ & = g_b ( q_{c,b} Y_1 +  \bar{q}_{a,b}Y_1, Y_2) (x)
\end{xalignat*} 
where $g_b$ is the metric  $\om(X, j_b Y)$. Since the Liouville form
$\mu_M$ is the
Riemannian volume for $g_b$, it comes that 
$$ \frac{ \delta(x)}{ \operatorname{det}^{\frac{1}{2}} [ -i \partial_{y^j} \partial_{y^k} \Phi
  (x,x,x)]_{j,k}} = \operatorname{det} ^{-\frac{1}{2}}[  q_{c,b} +
  \bar{q}_{a,b}] (x) $$
Thus to obtain the formula in theorem \ref{P1}, we have to show that
$$ \operatorname{det} ^{-\frac{1}{2}}[  q_{c,b} +  \bar{q}_{a,b}] =
\bi (j_a,j_b,j_c).$$

To see this, let us choose $j_b$ as the reference complex structure
and let us associate to $j_a$ and $j_c$ the bundle maps $\mu_a$ and
$\mu_c$ from $\La^{1,0}_{j_b}T^*M$ to $\La^{0,1}_{j_b}T^*M$ as in \eqref{eq:defmu}.
On one hand, we have by \eqref{eq:zeta} (since $\mu_b =0$)
$$ \bic (j_a,j_b,j_c) =
\operatorname{det}^{-1}  \left( \begin{array}{cc} \identite & \bar{\mu}_a \\
    \mu_c & \identite \end{array} \right) $$
On the other hand, $T^{0,1}_{j_a}M$ is the graph of $-\mu_a^t :
T^{0,1}_{j_b}M \rightarrow T^{1,0}_{j_b}M $. It follows that
$\bar{q}_{a,b}$ is the map
$$  \left( \begin{array}{cc} 0 & 0 \\
     \bar{\mu}_a^t & \identite \end{array} \right) : T^{1,0}_{j_b}M
     \oplus T^{0,1}_{j_b} M \rightarrow T^{1,0}_{j_b}M
     \oplus T^{0,1}_{j_b}M $$
Similarly, $q_{c,b}$ is the map 
$$   \left( \begin{array}{cc} \identite & \mu_c^t \\
     0 & 0 \end{array} \right) :  T^{1,0}_{j_b}M
     \oplus T^{0,1}_{j_b}M \rightarrow T^{1,0}_{j_b}M
     \oplus T^{0,1}_{j_b}M $$
The result follows.

\section{Half-form bundle and quantization} \label{sec:hfq}

\subsection{Preliminaries on half-form bundle} \label{sec:hfb}
Let $j$ be an almost-complex structure of $M$. Recall that a
half-form bundle of $(M,j)$ is a complex line bundle $\delta
\rightarrow M$ together with a line bundle isomorphism
 $$\varphi :
\delta^2 \rightarrow \La_j^{n,0} T^* M$$ 
which covers the identity of
$M$. Two half-form bundles $(\delta_a, \varphi_a)$ and $(\delta_b,
\varphi_b)$ are isomorphic if there exists a line bundle isomorphism
$\Psi : \delta_a \rightarrow \delta_b$ covering the identity and such that
 $$ \varphi_b \circ \Psi^2  = \varphi_a.$$ 
In the case where there exists a half-form bundle, there are $\# H^1(M,
\Z_2)$ isomorphism classes of half-form bundles.  
\scratch{rajouter un mot sur l'existence et donner une r{\'e}f{\'e}rence}

The existence and the choice up to isomorphism of a half-form bundle
over a symplectic manifold $(M, \om)$ is
in some sense independent of the almost complex structure, providing
it is compatible with $\om$ and positive. To see this we extend the previous
notion of half-form bundle isomorphisms to the collection $\demi$ consisting of the triples $(j, \delta, \varphi)$,
where $j$ is an almost-complex structure of $M$ compatible with $\om$
and positive, and $(\delta,
\varphi)$ is a half-form bundle for $(M,j)$. 

Let us define a half-form bundle morphism
$(j_a,\delta_a, \varphi_a) \rightarrow (j_b,\delta_b, \varphi_b)$ 
to be an isomorphism of line bundles $\Psi : \delta_a
\rightarrow \delta_b$ such that  
\begin{gather} \label{morphisme} 
 \varphi_b \circ \Psi^2 = \Psi_{j_a,j_b} \circ \varphi_a.
\end{gather} 
Here $\Psi_{j_a,j_b}$ is the morphism $\La_{j_a}^{n,0} T^*M
\rightarrow \La_{j_b}^{n,0} T^*M$ defined over $x \in M$ as in
definition \ref{def:Psi} with $E = T_x M$ and the complex structures $j_a(x)$ and $j_b (x)$.

The composition of a morphism $\Psi :(j_a,\delta_a, \varphi_a) \rightarrow
(j_b,\delta_b, \varphi_b)$ with a morphism $\Psi' :  (j_b,\delta_b, \varphi_b) \rightarrow
(j_c,\delta_c, \varphi_c)$ is defined as 
$$ \Psi' \compm \Psi :=  \bi (j_a,j_b,j_c) \Psi' \circ
\Psi  $$  
where the product $\circ$ on the
right-hand side is the usual composition of maps and the function $\bi
(j_a,j_b,j_c)$ is defined as in section \ref{section:prel}. Observe that $\compm$ is the product of symbol appearing in theorem \ref{P1}.

It is easily checked that $\demi$ with
this collection of  morphisms is a groupoid. The associativity
of $\compm$ follows from the
cocycle condition \eqref{eq:cocycle}. Equations \eqref{eq:ident} imply
that the identity $1_a$ of $\delta_a$ is the unit of $ (j_a,\delta_a,
\varphi_a)$, i.e.
$$ 1_a \compm \Psi = \Psi, \qquad \Psi'
\compm 1_a =  \Psi' ,$$
if $\Psi$ and $\Psi'$ are any morphisms $(j_b,\delta_b, \varphi_b)
\rightarrow (j_a,\delta_a, \varphi_a)$ and $(j_a,\delta_a, \varphi_a) \rightarrow
(j_b,\delta_b$,$ \varphi_b)$ respectively. Moreover, for any
$(j,\delta,\varphi) \in \demi$, define the Hermitian structure of
$\delta$ in such a way that $\varphi$ becomes an isomorphism of
Hermitian bundles. Then since $$\Psi_{j_a,j_b}^* =
\Psi_{j_b,j_a},$$ the adjoint $\Psi^*$ of any morphism $\Psi :
(j_a,\delta_a, \varphi_a) \rightarrow (j_b,\delta_b, \varphi_b)$ is a morphism $(j_b,\delta_b, \varphi_b) \rightarrow (j_a,\delta_a, \varphi_a) $ satisfying
\begin{gather} \label{inverse_adjoint} \Psi^* \compm \Psi
  = 1_a,  \qquad  \Psi \compm \Psi^* = 1_b. \end{gather} 
So $\Psi$ is invertible, with inverse $\Psi^*$. 

If $a$ and $b$ in $\demi$ are isomorphic, there exists exactly two
morphisms $a \rightarrow b$.  Observe also that given an almost complex
structure $j$, each isomorphism class of $\demi$ has a representative
whose almost complex structure is  $j$. So the existence of a half-form
bundle doesn't depend on the almost complex structure. And there are $\# H^1(M, \Z_2)$ isomorphism classes
in $\demi$ if it is not empty. 

\subsection{Quantization} \label{sec:Q}

Let us consider now the collection $\demint$ consisting of triples
$(j, \delta, \ph) \in \demi$ with
 an integrable complex structure $j$.
Given $a \in \demint$, let us denote by $\Quant_k (a)$ the Hilbert space of holomorphic sections
of $L^k \otimes \delta_a$. With our previous notations
 $$\Quant_k (a):= \Hilbert_k (j_a, \delta_a).$$  
Here the holomorphic and Hermitian structures of $\delta_a$ are such
that $\varphi_a :
\delta_a^2 \rightarrow \La_{j_a}^{n,0}T^* M$ is an isomorphism of
holomorphic Hermitian bundle. 

If $a$ and $b$ belongs to $\demint$, any
half-form bundle morphism $\Psi : \delta_a \rightarrow \delta_b$ is
the symbol of a unitary Fourier integral operator of $\Fourier ( (j_a,
\de_a), (j_b, \de_b))$
$$  \Um_k ( \Psi) : \Quant_k(a) \rightarrow \Quant_k (b), \quad k =1,2,... $$
Indeed if $S$ is a Fourier
integral operator with symbol $\Psi$, it follows from
\eqref{inverse_adjoint} and theorem \ref{P1} that $S^*S$ is a Toeplitz
operator with symbol 1. Hence $S
  (S^*S)^{-1/2}$ is a unitary Fourier integral operator with
  symbol $\Psi$. 

Contrary to the notations, $\Um (\Psi) = (\Um_k( \Psi))_k$ is not uniquely determined by
$\Psi$. It is unique modulo multiplication by a unitary Toeplitz
operator of symbol 1. So strictly speaking, $\Um ( \Psi )$  is an
equivalence class of 
Fourier integral operators. To avoid any confusion we will say that
two such operators are equal modulo $O(\hb)$. 

\begin{theo} \label{theo:fonct}
$\Um$ is functorial, that is if $\Psi''$ is the
  composition of the morphisms of half-form bundle $\Psi'$ and $\Psi$, then 
$$  \Um  ( \Psi'') = \Um (\Psi') \circ \Um
(\Psi) \text{ modulo }O(\hb).$$ 
Furthermore if $\Psi$ is a half-form bundle morphism $a
\rightarrow b$, the map sending the Toeplitz operator $T : \Quant (a)
\rightarrow \Quant (a) $ into $$(\Um(\Psi))^* T \Um (\Psi): \Quant (b)
\rightarrow \Quant (b) $$ preserves the normalized symbols modulo $O(\hb^2)$.
\end{theo}

First part is an immediate consequence of theorem \ref{P1} because the
composition of half-form bundle morphisms is the same as the
composition of symbols. Second part follows from theorem \ref{P2}, or
more directly from formula  \eqref{trans_symb}.

The group $\Group$ of connection-preserving Hermitian automorphisms of
$L$ acts on the quantum spaces as follows. First an automorphism $\Phi$ of
$\Group$ covers a symplectomorphism $\phi$ of $M$. Then $\Phi$ acts on
$\demint$  by sending $a=(j, \delta, \varphi)$
into $\Phi^*a := (\phi^*j,\phi^*\delta, \phi^*\varphi)$, where $ \phi^*\varphi$ is
defined in such a way that the diagram
\begin{gather*}
 \begin{CD} 
{\La_{j} ^{n,0} T^*_{\phi (x)} M } @>\phi^* >> {\La^{n,0}_{\phi^* j}
 T^* _x M} \\
 @A{\varphi}AA @AA{\phi^*{\varphi}}A \\
{\delta^2_{\phi(x)}} @>{(\phi^*)^2}>> {(\phi^* \delta)^2_x} \end{CD} \end{gather*}
commutes. Finally the operator  
 $$\Ci(M, L^k \otimes \delta_a) \rightarrow \Ci(M, L^k \otimes
 \phi^*\delta_{a}), \qquad s \rightarrow  ((\Phi^k)^{-1}
 \otimes \phi^*)\circ   s \circ \phi $$ 
restricts to a  unitary
 operator $\Quant_k (a) \rightarrow \Quant_k (\Phi^*a)$ that we denote by $\Phi^*$. 

Let us consider now $a \in \demi$, fixed until the end of this section. 
If $\Phi$ belongs to the identity component $\Group_o$ of $\Group$, then $a$
and $\Phi^* a$ are isomorphic half-form bundles. In this case, we
associate to any morphism $\Psi : \Phi^* a \rightarrow a$ the sequence
of operators
$$ \Vm_k ( \Psi, \Phi ) := \Um_k ( \Psi) \circ
\Phi^* : \Quant_k (a) \rightarrow \Quant_k (a)$$
As the operators $\Um( \Psi)$,  $ \Vm ( \Psi, \Phi )$ is uniquely defined up
to multiplication by a unitary Toeplitz operator with symbol 1. Denote
by $\Groupm$ the set of pairs $(\Psi, \Phi)$ where $\Phi
\in \Group_o$ and $\Psi$ is a half-form bundle morphism $\Phi^* a \rightarrow a$.

\begin{theo} \label{theo:repr}
For any half-form bundle $a \in \demint$, $\Groupm$ endowed
  with the product
$$ ( \Psi_1, \Phi_1).  ( \Psi_2, \Phi_2) :=   ( \Psi_2
\compm (\Phi_2^* \Psi_1) , \Phi_1 \circ \Phi_2) $$
is a central extension of $\Group_0$ by $\Z_2$. Furthermore $\Vm$
is a right-representation of $\Groupm$ up to $O(\hb)$ in the sense
that 
$$ \Vm ( \Psi_2, \Phi_2) \circ  \Vm ( \Psi_1,
\Phi_1) \equiv \Vm ( ( \Psi_1, \Phi_1).  ( \Psi_2,
\Phi_2)) \mod O(\hb).$$
\end{theo}

In the definition of the product of $\Groupm$, we used the following
action of $\Group$ on the half-form bundle morphisms. If
$\Phi$ is prequantization bundle automorphism of $L$ covering the
symplectomorphism $\phi$ and $\Psi$ is a morphism $a \rightarrow b$,
then $\Phi^* \Psi$ is the morphism $\Phi^* a \rightarrow \Phi^* b$
defined in such a way that the diagram 
\begin{gather*}
 \begin{CD} 
{\delta_a } @>\Psi >> {\delta_b} \\
 @V{\phi^*}VV @VV{\phi^*}V \\
{\phi^* \delta_a} @>{\Phi^*\Psi}>> {\phi^* \delta_b } \end{CD} \end{gather*}
commutes. One deduces easily from the relations
$$ \Phi_2^* (\Phi_1^* \Psi) = ( \Phi_1 \circ \Phi_2)^* \Psi, \qquad
\Phi^* ( \Psi_1 \compm \Psi_2 ) = (\Phi^*  \Psi_1)
\compm ( \Phi^* \Psi_2 ) $$
that $\Groupm$ is a group. Furthermore, one has 
$$  \Um ( \Phi^* \Psi ) = \Phi^* \circ \Um (\Psi )
 \circ (\Phi^*)^{-1} \mod O(\hb) $$
which implies the last part of the theorem:
\begin{xalignat*}{1}  
\Vm ( \Psi_2, \Phi_2) \circ  \Vm ( \Psi_1,
\Phi_1) = & \Um ( \Psi_2) \circ \Phi_2^* \circ
\Um ( \Psi_1) \circ \Phi_1^* \\
 = & \Um ( \Psi_2) \circ \Phi_2^* \circ
\Um ( \Psi_1) \circ ( \Phi_2^* )^{-1} \circ ( \Phi_1 \circ \Phi_2
)^* \\
 = &  \Um ( \Psi_2) \circ \Um (\Phi_2^* \Psi_1) \circ (\Phi_1
 \circ \Phi_2
)^* \mod O(\hb) \\
 = & \Um (\Psi_2 \compm  \Phi_2^* \Psi_1) \circ ( \Phi_1 \circ \Phi_2
)^*\mod O(\hb) \\
 = & \Vm ( (\Psi_1, \Phi_1 ).(\Psi_2,
\Phi_2))
\end{xalignat*}
by theorem \ref{theo:fonct}.

It is well-known that the Lie algebra of $\Group_0$ is $\Ci(M, \R)$,
the Lie bracket
being the Poisson bracket (cf. \eqref{eq:inf} for an explicit formula
for the exponential map). Let us associate to any $f \in \Ci(M)$ a Toeplitz operator $\Quant (f )$
of $\Quant (a)$ whose normalized symbol is $f$ modulo $O(\hb
^2)$. By theorem \ref{SP_norm}, we
obtain a Lie algebra representation up to $O(\hb)$ in the sense that 
$$ \tfrac{1}{i} \bigl[k \Quant_k (f), k\Quant_k (g) \bigr] = k\Quant_k (\{ f,g \} )
+ O(k^{-1}) $$
By exponentiating we recover the representation of theorem \ref{theo:repr}.

\begin{theo} \label{theo:exp}
For any $f \in \Ci (M, \R)$, we have
$$ \exp \bigl( it. k \Quant_k (f )  \bigr) = \Vm_k ( \Psi_t, \Phi_t) \mod O(\hb)$$
where $\Phi_t = \exp(tf)$ and $(\Psi_t)$ is the continuous family of half-form bundle morphisms
$\Psi_t: a \rightarrow \Phi_t^*a$ such that $\Psi_0$ is the identity
of $\delta_a$. 
\end{theo}
This last result will be proved in section \ref{sec:acG} (cf. remark after corollary
\ref{cor:Schrod}).

\subsection{Reformulation of the results of chapter \ref{FIO} }
 \label{reformulation}

Assume that $(M,\om)$ admits a unique half-form bundle
up to isomorphism. If this is not the case we can still apply what
follows by restricting to an open contractible set of $M$. 

Let us return to the quantum spaces $\Hilbert (a)$ defined from a
complex structure $j_a$ and a Hermitian holomorphic line bundle $K_a
\rightarrow M$. As in \cite{firstpart}, we introduce a half-form bundle
$(\delta_a , \varphi_a)$ and a holomorphic Hermitian line bundle
$L_{1,a}$ such that 
$$ K_a = \delta_a \otimes L_{1,a} .$$ 
For another pair $(j_b, K_b)$, introduce in the same way $(\delta_b,
\varphi_b)$ and $L_{1,b}$.
Then rewriting the formulas of chapter \ref{FIO} with these data, we get more transparent
results:
\begin{itemize}
\item
The formula for the commutators in theorem \ref{SP_norm} becomes
$$\tfrac{i}{\hb} (f *_a g - g *_a f) = \{ f, g \} - \hb \langle
\om_{1,a} , X_f \wedge X_g \rangle + O(\hb^2),$$
where $\frac{1}{i} \om_{1,a}$ is the curvature of $L_{1,a}$.
\item 
Denote by $\Mor_{a,b}$ the set of
half-form bundle morphisms $(j_a, \delta_a ,
\varphi_a) \rightarrow (j_b, \delta_b, \varphi_b)$, then  
$$ \Hom ( K_a, K_b) = \Mor_{a,b} \times_{\Z_2} \Hom(L_{1,a} , L_{1,b})$$
where we divided by $\Z_2$ to identify $(\Psi, \Psi_1)$ with $(-\Psi,
-\Psi_1)$. The composition of symbols in theorem \ref{P2} is then the
product of the composition of half-form bundle morphisms with the
usual composition. 
\item 
The symbol of a
unitary operator $S
\in U \Fourier (a,b)$ is of the form $[\Psi, \Psi_1 ]$ with $\Psi \in
\Mor_{a,b}$ and $\Psi_1$ a
unitary isomorphism $L_{1,a} \rightarrow L_{1,b}$. Furthermore the equivalence of the star-products $*_a$ and
$*_b$ induced by $S$ is up to second order 
$$ f_a \rightarrow f_a + \hb \langle \al_1, X_{f_a} \rangle + O(\hb^2) $$
where $\al_1$ is such that 
$ \nabla^{\Hom( L_{1,a}, L_{1,b} )} \Psi_1 = \tfrac{1}{i} \al_1 \otimes \Psi_1.$  
\end{itemize} 
This point of view will also be useful to prove theorem \ref{P2} in
the following section.

\section{Proof of theorem \ref{P2}} \label{preuve1} 

To prove the theorem, we consider the kernels of the Fourier integral
operators as Lagrangian sections and interpret $T_b S - S T_a$ as the
result of the action of a Toeplitz operator on a Lagrangian
section. The computation of the symbol is then a corollary of theorem
 3.4 in \cite{firstpart}. 

Let us regard $M^2$ as a symplectic manifold with symplectic form 
$$\om_{M^2} = \pi_l ^*
\om - \pi_r ^* \om,$$ 
where $\pi_l$ and $\pi_r$ are the projections onto the first and second
factor respectively. Then $L \boxtimes \bar{L}$  is a prequantization
bundle of $M^2$ with curvature $\frac{1}{i} \om_{M^2}$ and the
diagonal map $\De :
M \rightarrow M^2$ is a Lagrangian embedding. Furthermore $(j_b,
-j_a)$ is a complex structure of $M^2$ compatible with $\om_{M^2}$
and positive. Denote by
$\Hilbert ( b,-a)$ the associated Hilbert space
$$ \Hilbert (b,-a) = \bigl\{ \text{holomorphic sections of  } (L^k \boxtimes
  \bar{L}^{k}) \otimes (K_b \boxtimes \bar{K}_a) \rightarrow M^2 \bigr\} $$
Then the Fourier integral operators of $\Fourier (a,b)$ are defined in
  such a way that their kernel is a Lagrangian section of $\Hilbert (
  b,-a)$ associated to the diagonal. 

Let $S \in \Fourier (a,b)$ with kernel $S(.)$. Let $T_a$ and $T_b$ be
Toeplitz operators of $\Hilbert (a) $ and
 $\Hilbert (b)$ with normalized symbols $ f_a (., \hb)$
 and $f_b(.,\hb)$ respectively. Then it is easily checked that the kernel of $T_b S - S T_a$ is
 $T S(.)$ where $T$ is a Toeplitz operator of $\Hilbert ( b,-a)$ with
 normalized symbol
$$ g(x,y, \hb) = f_b (x, \hb) - f_a (y, \hb).$$
Assume now that $T_a$ and $T_b$ have the same principal symbol. Then
the principal symbol of $T$ vanishes along the diagonal and
consequently the principal symbol of $T S(.)$ vanishes. By applying 
theorem 3.4 of \cite{firstpart}, we obtain the principal symbol of $k^{-1}T S(.)$ which
corresponds to the principal symbol of $k^{-1} (  T_b S
- S T_a)$. 

Let us use the half-form bundles as in chapter \ref{reformulation}.  
The symbol of $S$ as a Fourier integral operator of $\Fourier (a,b)$
is a class $[\Psi,
\Psi_1 ]$ where  $\Psi$ is a half-form bundles morphism $\delta_a \rightarrow
\delta_b$ and $\Psi_1$ a bundle morphism $L_{1,a} \rightarrow
L_{1,b}$. We have to show that the symbol of $k^{-1} (  T_b S
- S T_a) $  is $ [\Psi,\Psi_1' ]$ with
\begin{gather} \label{but} 
 \Psi'_1 = (f_{1,b} - f_{1,a}) \Psi_1 + \tfrac{1}{i}
\nabla^{\Hom ( L_{1,a}, L_{1,b} )} _X \Psi_1 \end{gather} 
where $X$ is the Hamiltonian vector field of $f$ and 
$f_{1,a}$, $f_{1,b}$ are the subprincipal symbols of $T_a$ and
$T_b$ respectively.

To the morphisms $\Psi$ and $\Psi_1$ correspond two sections 
$$ \tilde{\Psi} \in \Ci (M, \delta_b \otimes
\overline{\delta}_a ), \quad 
\tilde{\Psi}_1 \in \Ci (M,  L_{1,b} \otimes \overline{L}_{1,a}).$$
The principal symbol of the Lagrangian section $S(.)$ is $ \tilde{\Psi}
\otimes \tilde{\Psi}_1$. The restriction to the diagonal of the
Hamiltonian vector field of the principal symbol of $T$ is $\Delta_*
X$. Then it follows from theorem 3.4 of \cite{firstpart} that the principal symbol of $k^{-1} T. S(.)$ is
$$  (f_{1,b} - f_{1,a}) \tilde{\Psi}
\otimes \tilde{\Psi}_1 + \tfrac{1}{i}  (D_X^{\De} \tilde{\Psi})
\otimes \tilde{\Psi}_1 + \tfrac{1}{i} \tilde{\Psi}
\otimes \nabla_X^{ L_{1,b} \otimes \overline{L}_{1,a}} \tilde{\Psi}_1$$
It remains to explain how is defined  the section $D_X^{\De}
\tilde{\Psi}$ of $\delta_b \otimes
\overline{\delta}_a$ and to
prove that it vanishes. This will imply
\eqref{but}.

Consider the isomorphism 
$$ \xi: \La_{j_b}^{n,0} T^*M \boxtimes \overline{\La_{j_a}^{n,0} T^* M}
\rightarrow \La ^{2n, 0}_{j_b, -j_a}T^*  M^2, \quad \be \boxtimes \overline{\al} \rightarrow
\pi_l^* \be \wedge \pi_r^* \overline{\al} $$
$(\delta_b \boxtimes \bar{\delta}_a, \xi \circ (\varphi_b \boxtimes
\bar{\varphi_a}))$ is a half-form bundle of $M^2$ for the complex
structure $(j_b, -j_a)$.  
Then $$ \delta_b \otimes \bar{\delta}_a = \De^* (\delta_b \boxtimes
\bar{\delta}_a)$$
is a square root of  $ \La^{2n} T^* M \otimes \C $ through the map
$$ \varphi_{\delta} : \delta_b^2 \otimes \bar{\delta}_a^2 \rightarrow 
 \La^{2n} M \otimes \C, \quad u_b \otimes u_a \rightarrow 
\De^* (\xi( \varphi_b (u_b) \boxtimes \bar{\varphi_a}(u_a)  )).$$
and $D_X^{\De} \tilde{\Psi}$ is defined in such a way that 
$$ \varphi_{\delta} ( \tilde{\Psi} \otimes D_X^{\De} \tilde{\Psi} ) =
\tfrac{1}{2} \Lie_X. 
\varphi_{\delta} (\tilde{\Psi}^{\otimes 2} )  $$
Then $D_X^{\De} \tilde{\Psi} = 0$ follows from the following lemma and
Liouville theorem.
\begin{lemme} 
$ \varphi_{\delta} (
\tilde{\Psi}^{\otimes 2} ) =i ^{n(n-2)} \om^n/ n!$ \end{lemme} 

\begin{proof} 
Denote by 
$\tilde{\Psi}_{j_a,j_b}$ the section of  
$ 
\La_{j_b}^{n,0} M \otimes \overline{\La_{j_a}^{n,0} M}
\rightarrow  M$
associated to $\Psi_{j_a,j_b}$. Since $\Psi$ is a half-form bundle
morphism, we have  
$$ (\varphi_a \otimes \bar{\varphi}_b) (\tilde{\Psi}^{\otimes 2}) =
\tilde{\Psi}_{j_a,j_b}$$
Introduce a unitary section $\al $ of $\La_{j_a}^{n,0} M$. We have   
$ \tilde{\Psi}_{j_a,j_b} = \Psi_{j_a,j_b} (\al ) \otimes \bar{ \al}$. 
Consequently 
\begin{xalignat*}{2}
  \varphi_{\delta} (
\tilde{\Psi}^{\otimes 2} )  = &
\De^* ( \pi_l^* \Psi_{j_a,j_b} (\al ) \wedge \pi_r^* \bar{ \al})) \\
= & \Psi_{j_a,j_b} (\al ) \wedge \bar{\al} \\
= & \al \wedge \bar{\al} 
\intertext{by definition of $\Psi_{j_a,j_b}$} 
=&  i ^{n(n-2)} \om^n/ n!
\end{xalignat*}
because $\al$ is unitary. \end{proof}

\section{Geometric interpretation} \label{sec:geo} 

Consider as previously a symplectic manifold $(M, \om)$ with a
prequantization bundle $(L, \nabla)$. The space $\Comp $ of almost
complex structures of $M$ compatible with $\om$ and positive may be
regarded as the space of sections of a fiber bundle over $M$, which
turns it into an infinite dimensional manifold. 
Let us fix a isomorphism class $D$ of half-form bundles and choose
for any $j \in \Comp$ a half-form bundle of $(M,j)$ which represents $D$ and depends "smoothly" on $j$. One
way to do that is first to choose $(j_0,\delta_{0},
\varphi_{0})$ representing $D$ and then to set
$$\delta_{j}:=\delta_{0}, \qquad \varphi_{j} := \Psi_{j_0, j} \circ
\varphi_{0}, \qquad \forall j \in \Comp.$$  
Let $\Preq_k \rightarrow \Comp$ be the bundle of
{\em prequantum} spaces, whose fiber
at $j$ is the space of smooth sec\-tions of $L^k \otimes
\delta_j$. Let us consider now a submanifold $\Compint$ of $\Comp$ which
contains only integrable complex structures. Assume that the family of Hilbert
spaces 
$$ \Quant_{k,j} := \{ \text{holomorphic sections of } L^k \otimes \delta_j
\}, \qquad j \in \Compint$$
defines a smooth subbundle $\Quant_k \rightarrow \Compint$ of $ \Preq_k
\rightarrow \Compint$, when $k$ is sufficiently large. This assumption is satisfied as soon as the
dimension of $ \Quant_{k,j}$ is constant when $j$ runs over
$\Compint$. This follows from Fredholm theory because $ \Quant_{k,j}$ is
the kernel of the holomorphic Laplacian, an elliptic second order
differential operator whose coefficient depend smoothly on the complex
structure. Furthermore as noticed by Foth and Uribe \cite{FoUr}, for
any complex structure $j_0$, there exists an integer $N$ such that the
dimension of $\Quant_{k,j}$ is constant when $j$ describe a $C^2$
neighborhood of $j_0$ and $k$ is larger than $N$. The $C^2$-topology
is involved here to control the curvature term in the Bochner-Kodaira
identity and deduce a uniform vanishing theorem. Then the dimension of
$\Quant_{k,j}$  is given by the Riemann-Roch theorem.

Before we continue, let us note that $\Preq_k$ and $\Quant_k$ depend only on the
isomorphism class $D$, providing we regard them as the orbifold bundles
$\Preq_k / \Z_2$ and $\Quant_k/ \Z_2$, where $\Z_2$ acts trivially on the
base $\Comp$ and by $\pm \identite$ on the fibers. Indeed, let us
consider another smooth family $(\tilde{\delta}_{j},\tilde{\varphi}_j)_{j\in \Comp}$,
obtained as above by choosing a half-form bundle $(\tilde{j}_0,\tilde{\delta}_{0},
\tilde{\varphi}_{0})$ representing $D$ and denote by $\tilde{\mathcal{P}}^{\operatorname{m}}_k$ the associated
bundle of prequantum spaces. 
Then there exists exactly two continuous families 
$$ \bigl( \Psi_j : (\delta_j, \varphi_j)
\rightarrow (\tilde{\delta}_{j}, \tilde{\varphi}_j); \;j\in \Comp
\bigr)$$ of half-form bundle morphisms. These families
induce isomorphisms $\Preq_k \rightarrow \tilde{{\mathcal{P}}}^{\operatorname{m}}_k$ and $\Quant_k \rightarrow
    \tilde{{\mathcal{Q}}}^{\operatorname{m}}_k$, which are unique up to the
$\Z_2$-action. All the constructions which follow only depend on $D$
in this sense.

First we define a connection on $\Preq_k$. Given a tangent vector $\mu$
of $\Comp$ at $j_0$, let us introduce a curve $j_t$ of $\Comp$ tangent to $\mu$
at $t=0$ and consider the continuous family $(\Psi_t)$ of half-form
bundle morphism $(\delta_{j_0}, \varphi_{j_0}) \rightarrow (\delta_{j_t},
\varphi_{j_t})$ such that $\Psi_0$ is the identity of $\delta_{j_0}$. Then we define the covariant
derivative of a section $\Phi$ of $\Preq_k$ with respect to $\mu$
to be 
$$ \nabla_{\mu}^{\Preq_k} \Phi \ (j_0) :=  \frac{d}{dt}
\Bigr|_{t=0} \; \Psi_t^{-1} . \Phi (j_t) $$
where the derivative is in the $t$-independent space
$\Preq_{k,j_0}$. The connection on $\Quant_k$ is then defined as 
$$ \nabla^{\Quant_k} := \Pi_k \circ \nabla^{\Preq_k} $$
where $\Pi_k$ is the section of $\End (\Preq_k )$ which at $j$ is the
orthogonal projector onto $ \Quant_{k,j}$.

\begin{theo} \label{theo:courbureQ}
For any $k$, the connection $ \nabla^{\Quant_k}$ is compatible with the
  Hermitian structure. Furthermore,
\begin{itemize}
\item 
  For any $j \in \Compint$ and $\eta,
  \mu \in T_j \Compint$, the sequence of curvature 
 $$R^{\Quant_k} ( \eta,\mu) : \Quant_{k,j} \rightarrow \Quant_{k,j},
 \quad k=1,2,...$$  
 is a Toeplitz operator
  whose principal symbol vanishes. 
\item
  For any curve $\ga$ of $\Compint$ with endpoints $j_a$ and $j_b$,
  the sequence of parallel transport $\ga$ in $\Quant_k$  is a
  unitary Fourier integral operator 
 $$ \Quant_{k,j_a} \rightarrow
  \Quant_{k,j_b}, \quad k=1,2,...$$ 
of $\Fourier( (j_a, \delta_{j_a}),( j_b, \delta_{j_b} ) )$. 
  Its principal symbol is the half-form bundle
  morphism $\delta_{j_a} \rightarrow \delta_{j_b}$ obtained by extending
  continuously the identity of $\delta_{j_a}$ in half-form bundle
  morphisms $ \delta_{j_a} \rightarrow \delta_{\ga (t)}$.
\end{itemize} \end{theo} 

The proof is postponed to section \ref{sec:preuve}. Let us compute the
curvature of $\Preq_k$. 
Given an almost complex-structure
$j_0 \in \Comp$, we can represent any $j \in \Comp$ as a section 
$$\mu \in \Ci (M,
\Hom( \La_{j_0}^{1,0}T^*M , \La_{j_0}^{0,1}T^*M))$$ 
such that the graph of $\mu(x)$ is $\La^{1,0}_j T^*_x M$ for any $x \in
M$. In this way, we identify the tangent space to $\Comp$ at $j_0$ with 
$$ T_{j_0} \Comp \simeq\bigl\{ \mu \in \Ci(M, \Hom( \La_{j_0}^{1,0}T^*M , \La_{j_0}^{0,1}T^*M)) ;
\ \om ( \mu^t (.),.) + \om (.,\ \mu^t(.) ) = 0 \bigr\} $$  
and $\Comp$ becomes a neighborhood of the
zero section of $ T_{j_0} \Comp$.

\begin{theo} \label{theo:courbureP}
The connection $ \nabla^{\Preq_k}$ is compatible with the
  Hermitian structure. Its curvature at $\eta,
  \mu \in T_j \Comp$ is given by
$$ R^{\Preq_k} ( \eta, \mu)\ \Phi = \frac{1}{2} \trace(
  \eta.\bar{\mu} - \mu. \bar{\eta} )\ \Phi, \qquad  \Phi \in \Preq_{k,j} $$
\end{theo} 

It is interesting to compare the previous theorems with the results
of Foth and Uribe \cite{FoUr}. The curvature of $\Quant_k$ is the sum of two terms  which cancel each other
at first order. The first term is the curvature of $\Preq_k$ and the second
one is a commutator (cf. lemma \ref{lem:imp}). In the case considered
by Foth and Uribe, the prequantum spaces are defined without half-form
and consequently don't depend on the complex structure. Then the bundle
 $\Compint \times \Ci (M, L^k)$ is endowed with the trivial connection,
 and composing with the Szeg{\"o} projector, we obtain a connection on
 the quantum space bundle. Its curvature equals a commutator
 (cf. lemma 2.1 of \cite{FoUr}), which is
 essentially the same as in our situation, and isn't canceled by the curvature of the
 prequantum bundle, flat in this case.

\begin{proof}[Proof of theorem \ref{theo:courbureP}]
Let $j_0$ be a fixed almost-complex structure and let us identify
$\Comp$ with an open convex set $O$ of $T_{j_0} \Comp$ as
previously.  Let us compute
the connection in the trivialization 
$$ \Preq_k \simeq O \times \Preq_{k,j_0} $$
induced by the continuous family of half-form bundle isomorphisms
$(\delta_{j_0} , \varphi_{j_0}) \rightarrow (\delta_{j} ,
\varphi_{j})$ extending the identity of $\delta_{j_0}$.

Let $\mu(t)$ be a curve of $O$
covered by a section $\Phi (t)$. By lemma \ref{lem:the_Lemme}, the
continuous curve of half-form bundle morphisms $\Psi_{t} :
\delta_{\mu(0)}  \rightarrow \delta_{\mu(t)}$ is in the previous
trivialization the multiplication by the continuous square root of 
$$  t \rightarrow \det \left( \begin{array}{cc} \identite & \bar{\mu}(0) \\
    \mu(0) & \identite \end{array} \right) 
. \operatorname{det}^{-1}  \left( \begin{array}{cc} \identite & \bar{\mu}(0) \\
    \mu(t) & \identite \end{array} \right) $$ 
equal to 1 at $t=0$. Then we have  
\begin{xalignat*}{2}
 \nabla_{\dot{\mu}(0)}^{\Preq_k} \Phi (0) & =  \frac{d}{dt} \Bigr|_{t=0}
 \Psi_t^{-1} . \Phi (t) \\ 
& = - \frac{1}{2} \trace \Bigl( \dot{\mu} (0) \bar{\mu} (0) \  \bigl( \identite -  \mu
 (0) \bar{\mu} (0) \bigr)^{-1} \Bigr) \Phi (0) + \dot{\Phi} (0)
\end{xalignat*} 
Thus we have $ \nabla^{\Preq_k} =  d + \al$ with 
$$ \langle \al, \dot{\mu}  \rangle. \Phi = - \frac{1}{2} \trace \Bigl( \dot{\mu}  \bar{\mu}  \  \bigl( \identite -  \mu
  \bar{\mu}  \bigr)^{-1} \Bigr) \Phi, \qquad \forall \ \dot{\mu} \in
  T_{\mu} O, \ \Phi \in \Preq_{k,j_0} $$
Finally it is easy to compute the curvature at the origin of $O$,
  where $\al$ vanishes which leads to the formula of the theorem. 

To check the compatibility of the connection with the scalar product,
observe that our trivialization doesn't preserve the scalar
product. Actually since $\Psi_{j_0,j}^* \circ \Psi_{j_0,j} = \bic
(j_0,j,j_0)$, we have 
$$
 ( \Phi, \Phi' )_{\Preq_{k,j}} = \int_M \frac{1}{\bi (j_0,j,j_0) (x) }
 \bigl( \Phi (x) , \Phi'(x) \bigr)_{L^{k} \otimes \delta_{j_0}} \mu_M (x)
$$
Then using again that the connection form $\al$ vanishes at the origin and that $\bic
(j_0,j,j_0) = 1 + O (|j - j_0|^2)$  we deduce that for every section
$\Phi, \Phi'$ of $ O \times \Preq_{k,j_0}$, 
$$  d ( \Phi, \Phi') = ( \nabla^{\Preq_k} \Phi, \Phi' ) + ( \Phi ,
\nabla^{\Preq_k} \Phi')$$
at the origin. \end{proof} 

\begin{rem}   
It is immediate to deduce the first part of theorem
\ref{theo:courbureQ}. Since $\Pi$ is self-adjoint, if $\Phi$ and
$\Phi'$ are section of $\Quant_k \rightarrow \Compint$, then 
\begin{xalignat*}{2}
  ( \nabla^{\Quant_k} \Phi, \Phi' ) + ( \Phi ,
\nabla^{\Quant_k} \Phi') =&  ( \Pi \nabla^{\Preq_k} \Phi, \Phi' ) + ( \Phi ,
\Pi \nabla^{\Preq_k} \Phi') \\
= &  ( \nabla^{\Preq_k} \Phi, \Phi' ) + ( \Phi ,
 \nabla^{\Preq_k} \Phi') \\
= &  d ( \Phi, \Phi') \end{xalignat*} 
which proves that $\nabla^{\Quant_k}$ is Hermitian. \qed
\end{rem}

\section{Action of the prequantization bundle automorphisms} \label{sec:acG}

Adapting the constructions of section \ref{sec:hfq}, we define an action of
the identity component $\Group_0$ of the group of prequantization bundle
automorphisms of $(L, \nabla)$ on $\Preq_k$ and $\Quant_k$. For any equivariant vector bundle equipped with an invariant
connection, one defines a moment (cf. definition 7.5 in
\cite{BeGeVe}). This notion makes sense in our infinite dimensional
setting. In the first part of this section, we prove the moment of a function $f$ in the Poisson
Lie algebra $\Ci (M)$ is a Toeplitz operator. From this we compute the
solution of Schr{\"o}dinger equation in terms of parallel
transport. This last result was obtained in \cite{FoUr} in the case
without metaplectic correction. This enables us to deduce that the quantum
propagator is a Fourier integral operator from the fact that parallel
transport is such an operator. Next we compute the commutator of
Toeplitz operators in terms of the curvature of the quantum space
bundle. This prove our estimate of the curvature is sharp. Finally we
explain how the same ideas apply in the case without metaplectic
correction.

\subsection{The infinitesimal action of $\Group_o$ on the bundles $\Preq$ and $\Quant$}
Let us start with the definition of the action of $\Group_0$ on $\Preq_k$. An
automorphism $\Phi \in \Group_o$ covering the
symplec\-tomorphism $\phi$ acts on the base $\Comp$ by sending $j$
into $\phi^* j$.  Let us lift this action to $\Preq_k$.  Given  $j \in \Comp$, there are exactly two bundle maps
$\phi^* : \delta_j \rightarrow \delta_{\phi^* j}$ covering $\phi$ such
that the diagram 
\begin{gather*}
 \begin{CD} 
{\La_{j} ^{n,0} T^*M} @>\phi^* >> {\La^{n,0}_{\phi^* j}
 T^*M} \\
 @A{\varphi_j}AA @AA{\varphi_{\phi^* j}}A \\
{\delta^2_{j}} @>{(\phi^*)^2}>> {\delta^2_{\phi^*j}} \end{CD} \end{gather*}
commutes. Then the pull-back by $$\Phi^k \otimes (\phi^*)^{-1} : L^k
 \otimes \delta_j \rightarrow L^k \otimes \delta_{\phi^*j}$$ is a
 linear map $\Preq_{k,j} \rightarrow \Preq_{k,\Phi^*j}$. Choosing the bundle
 map $\delta_j \rightarrow \delta_{\Phi^* j}$ in such a way that it
 depends continuously on $j$, we obtain the action of $\Phi$
 on $\Preq_k$.
Since this action is only defined up to multiplication by $-1$, we obtain merely a $\Group_o$-action  on the orbibundle $(\Preq_k, \Z_2)$.

Given a function $f \in \Ci (M)$, let us define the operator  
$$ \Op_k (f) := f + \frac{1}{ik} ( \nabla^{L^k}_{X}
\otimes \identite + \identite \otimes D_X ): \Preq _{k,j} \rightarrow
\Preq_{k,j},  \qquad j \in \Comp$$ 
where $X$ is the Hamiltonian vector field of $f$ and  $D_X$ is
the first order differential operator of $\Ci (M, \delta_j)$ such that 
\begin{gather} \label{eq:defDX}
  p_{j} \Lie_X   \varphi_j ( \be ^2) = 2 \varphi_j \bigl( \be
\otimes (D_X \be) \bigr), \qquad \forall \ \be \in \Ci (M, \delta_{j})
\end{gather} 
with $p_{j}$ the projection of $\La^nT^*M \otimes \C$ onto $\La_{j}^{n,0}T^*M $ with
kernel $ \La_{j}^{n-1,1} T^*M \oplus ...\oplus \La^{0,n}_{j}T^*M $. 

Recall that the Lie algebra of $\Group_o$ may be viewed as $\Ci(M, \R)$. Given $f \in
\Ci(M,\R)$, $\Phi_t := \exp (tf)$ is the automorphism of $L$ which covers the Hamiltonian flow
$\phi_t$ of $f$ at time $t$ and is given by
\begin{gather} \label{eq:inf}
 \Phi_t (\xi)  = e^{i t f(x) } {\mathcal{T}}_{t} .\xi,
\qquad \text{ if } \xi \in L_x \end{gather}  
where ${\mathcal{T}}_{t}$ is the parallel transport from $L_x$ to
$L_{\phi_t(x)}$ along the Hamiltonian flow.

\begin{theo} \label{theo:infinit}
Let $f \in \Ci (M,\R)$ and denote by 
$$U_t: \Preq_k \rightarrow \Preq_k$$ 
the action of $\exp (tf)$ on $\Preq$.
Let $j_0 \in \Comp$ and $j: \R \rightarrow \Comp$ be the curve $j_t
= \phi_t^* j_0 $. For any $s _0 \in \Preq_{k,j_0}$ the
section $s$ of $j^*\Preq_k$ defined by $s_t = U_t.s_0$ satisfies
$$ \frac{1}{ik}\nabla_{\partial_t} ^{j^* \Preq_k} \ s = \Op_k (f)\ s . $$
\end{theo}

The action of the symplectomorphism group on $\Comp$ preserves the subspace of integrable almost complex
structures. Assume that $j_0 \in \Compint$ and that the submanifold $\Compint$ is invariant
under the action of $\exp (tf)$. Then the operator $U_t$ restricts to
an operator $\Quant_{k,j_0} \rightarrow \Quant_{k,j_t}$ and the section
$s$ of $ j^*\Quant_k$ defined as above satisfies
$$ \frac{1}{ik}\nabla_{\partial_t} ^{j^* \Quant_k} \ s = \Quant_k (f) \ s
$$
where $ \Quant_k (f)$ is defined by 
$$ \Quant_k (f) := \Pi_{k,j}  \Op_k(f) : \Quant_{k,j} \rightarrow \Quant_{k,j} .$$
We proved in \cite{firstpart} (cf. theorem 1.5) the the sequence
$(\Quant_k (f))_k$ is a Toeplitz operator whose normalized symbol is $f + O( \hb^2)$.

\begin{proof}
First we may assume that $s_0 = \al \otimes \be$ with $\al \in \Ci
(M, L^k)$ and $\be \in \Ci(M,\delta_{j_0})$. Furthermore since $U_t \circ U_s = U_{t+s}$, it
is sufficient to prove the result at $t =0$.   
Let us write $$\al_t
= \Phi^*_t  \al, \qquad \be_t = \phi^*_t \be \in
\Ci(M, \delta_{j_t}).$$ Then $U_t. (\al \otimes \be) = \al_t \otimes \be_t$ and
consequently  
$$ \nabla_{\partial_t} ^{j^* \Preq_k} \ s =
\dot{\al}_t \otimes \be_t + \al_t \otimes \dot{\be}_t $$
where the derivative $\dot{\al}_t$ is in the $t$-independent vector
space $\Ci(M, L^k)$ and 
\begin{gather} \label{eqeq} 
 \dot{\be}_0 := \frac{d}{dt}\Bigr|_{t =0} \Psi_t^{-1} \Phi^*_t \be  
\end{gather}
with $\Psi_t : \delta_{j_0} \rightarrow \delta_{j_t}$ the
continuous family of half-form bundle morphisms such that $\Psi_0$ is
the identity of $\delta_{j_0}$. It is a classical result that 
$$ \dot{\al}_t = ik \Bigl( f +  \frac{1}{ik} \nabla^{L^k}_{X_f} \Bigr) \al_t$$
So it remains to prove that 
$$ \dot{\be}_0 = D_X  \be $$
Denote by $\tilde{D}$ the map sending $\be \in \delta_{j_0}$
into $\dot{\be}_0$ defined in \eqref{eqeq}. $\tilde{D}$ is a first order differential
operator. We have to prove that 
$$ p_{j_0} \Lie_X \varphi_{j_0} ( \beta^2) = 2 \varphi_{j_0} ( \beta \otimes
\tilde{D} \beta)$$
We have 
\begin{xalignat*}{2} 
 p_{j_0} \Lie_X \varphi_{j_0}  ( \beta^2) = &  p_{j_0} \frac{d}{dt}
 \Bigr|_{t=0}  \phi_t^*
\varphi_{j_0}  ( \beta^2) \\
 = & p_{j_0} \frac{d}{dt} \Bigr|_{t=0} \varphi_{j_t}
 (( \phi_t^*\beta)^2)  \\
=& \frac{d}{dt}  \Bigr|_{t=0} p_{j_0} \varphi_{j_t}  (( \phi_t^*\beta)^2)
\end{xalignat*} 
Now it follows from the definition of $\Psi_{j_0,j_t}$ that 
$$ p_{j_0} \ga = \Psi_{j_0,j_t}^{-1} \ga, \qquad \forall \ \ga \in
\La^{n,0}_{j_t} T^*M $$
Since $\Psi_t $ is a half-form bundle morphism, we obtain 
\begin{xalignat*}{2} 
 p_{j_0} \Lie_X \varphi_{j_0}  ( \beta^2) 
= & \frac{d}{dt} \Bigr|_{t=0} \Psi_{j_0,j_t}^{-1}
 \varphi_{j_t} ( (\phi_t^*\beta)^2)  \\
= & \frac{d}{dt} \Bigr|_{t=0}
 \varphi_{j_0} ( ( \Psi_t^{-1} \phi_t^*\beta)^2) \\
= &  \varphi_{j_0} \frac{d}{dt} \Bigr|_{t=0} ( ( \Psi_t^{-1} \phi_t^*\beta)^2)
  \\
= & 2 \varphi_{j_0} ( \beta \otimes
\tilde{D} \beta)
\end{xalignat*} 
which was to be proved. 
\end{proof}

\subsection{Schr{\"o}dinger equation} 

As a corollary of theorem \ref{theo:infinit}, we obtain the relation between the parallel transport
in $\Preq_k$, the action of $\Group_0$ and the Schr{\"o}dinger equation with
Hamiltonian $\Op_k (f)$.
\begin{cor} \label{cor:Prop}
Let ${\mathcal{T}}_t^{\Preq_k} : \Preq_{k,j_0} \rightarrow
\Preq_{k,j_t}$ be the parallel transport along the curve $t \rightarrow
j_t$. Then the family of operators $$P_t :=  ({\mathcal{T}}^{\Preq_k}_t)^{-1}\circ U_t : \Preq_{k,j_0} \rightarrow
\Preq_{k,j_0}$$ satisfies 
$$ \frac{1}{ik} \frac{d}{dt} P_t s_0 = \Op_k (f) \ P_t s_0 $$
for any $s_0 \in \Preq_{k,j_0}$.
\end{cor}

\begin{proof} Since ${\mathcal{T}}_{t}^{\Preq_k}$ is parallel transport,
\begin{xalignat*}{2} 
\frac{1}{ik} \frac{d}{dt} P_t s_0 = & \frac{1}{ik} 
({\mathcal{T}}_{t}^{\Preq_k})^{-1} \nabla_{\partial_t} ^{j^* \Preq_k}  U_t
s_0 \\ = & ({\mathcal{T}}_{t}^{\Preq_k})^{-1} \Op_k (f) U_t s_0
\end{xalignat*} 
by theorem \ref{theo:infinit}.
Furthermore for any function $g \in \Ci(M)$, we have
$$ U_t \Op_k ( g ) = \Op_k ( \Phi_t^* g ) U_t .$$
So $U_t$ and
$\Op_k( f)$ commutes, because $f$ is preserved by its Hamiltonian flow. 
Consequently 
\begin{gather} \label{eq:clef}
 \frac{1}{ik} \frac{d}{dt} P_t s_0 = P_t \Op_k (f) s_0 . \end{gather} 
To conclude, we prove that  $P_t$ and $\Op_k (f)$
commute. We have
\begin{xalignat*}{2} 
 \frac{d}{dt} \bigl( P_t \Op_k (f) (P_t)^{-1} \bigr) = & \dot{P}_t \Op_k
(f) (P_t)^{-1} -  P_t \Op_k (f) (P_t)^{-1} \dot{P}_t  (P_t)^{-1} \\
=& ik  \bigl (P_t \Op_k
^{2} (f) (P_t)^{-1} - P_t \Op_k (f) (P_t)^{-1} P_t \Op_k (f)
(P_t)^{-1}\bigr)  \\ =&  0 \end{xalignat*}
because of \eqref{eq:clef}.
\end{proof} 

Let us assume again that $j_0 \in \Compint$ and that $\Compint$ is
preserved by the action of $\exp (tf)$. Then arguing as in the
previous proof, we can deduce the similar result for the bundle
$\Quant_k$. 

\begin{cor} \label{cor:Schrod}
Let ${\mathcal{T}}^{\Quant_k}_t : \Quant_{k,j_0} \rightarrow
\Quant_{k,j_t}$ be the parallel transport in $\Quant_k$ along the
curve $t \rightarrow j_t$. Then we have 
$$   \frac{1}{ik} \frac{d}{dt} {P}^{\Quant_k}_t s_0 = \Quant_k (f)
{P}^{\Quant_k}_t  s_0 , \qquad \forall s_0 \in \Quant_{k,j_0}$$
where ${P}^{\Quant_k}_t$ is the operator
$({\mathcal{T}}^{\Quant_k}_t)^{-1} \circ  U_t : \Quant_{k,j_0} \rightarrow
\Quant_{k,j_0}$. 
\end{cor}

Recall that $\Quant (f)$ is a Toeplitz operator whose normalized
symbol is $f + O( \hb^2)$. Then theorem \ref{theo:exp} follows from the fact that the
parallel transport in $\Quant$ is a unitary Fourier integral operator
(cf. theorem \ref{theo:courbureQ}).

\subsection{Commutators and curvature} 

In the next theorem, we 
compute the commutator of $\Op_k (f)$ and $\Op_k (g)$ (resp. $\Quant_k (f)$
and $\Quant_k (g)$) in terms of the
curvature of $\Preq_k \rightarrow \Comp$ (resp. $\Quant_k \rightarrow
\Comp$). 

\begin{theo} \label{theo:courbcom}
Let $f$ and $g$ be two functions of $\Ci(M,\R)$. Then
$$ ik [ \Op_k (f) , \Op_k (g) ] = \Op_k ( \{ f,g \} ) + (ik)^{-1} R^{\Preq_k}(\eta, \mu)$$
where $\eta$ and $\mu$ are the vector fields of $\Comp$ corresponding
to the infinitesimal action of $ f$ and $g$ on $\Comp$. Furthermore, 
$$ ik [ \Quant_k (f) , \Quant_k (g) ] = \Quant_k ( \{ f,g \} ) + (ik)^{-1} R^{\Quant_k}(\eta, \mu)$$
when $\eta$ and $\mu$ are tangent to $\Compint$.
\end{theo} 

Since $\Quant (f)$ is a Toeplitz operator with normalized symbol $f +
O(\hb^2)$, it follows from  theorem \ref{SP_norm} that
\begin{gather} \label{eq:com}
ik [ \Quant_k (f) , \Quant_k (g) ] = \Quant_k ( \{ f,g \} ) + O(
k^{-2}) \end{gather} 
This is consistent with the second equation of the previous theorem
and  the fact that
$R^{\Quant_k}(\eta, \mu)$ is $O(k^{-1})$. Moreover, on can prove in
this way that
$R^{\Quant_k}$ can't be $O(k^{-2})$ except for particular sub manifolds
$\Compint$. Indeed given a complex structure
$j$, there is a star-product $*_j$ such that for any functions $f$ and
$g$, 
$$ \Quant_k (f)_j \circ \Quant_k (g)_j \equiv \Quant_k (h(.,k))_j +  O(k^{-\infty})$$
where $h(.,k)$ has an asymptotic expansion  $h_0 + k^{-1} h_1 +..$
whose coefficients satisfy $f *_j g = \sum \hbar^l h_l$. One can prove
that $*_j$ is a Vey star-product, i.e. the bidifferential operators
defining $*_j$ have the same principal symbol than the bidifferential
operators defining the Moyal-Weyl star-product. Hence
$$ 
i \hb^{-1} \bigl( f *_j g - g  *_j f \bigr) =  \{ f,g \} + \hb^2 A (f,g) + O( \hb^3)
$$
where $A$ is a non-vanishing bidifferential operator. So if $\eta$ and
$\mu$ are the infinitesimal actions of $f$ and $g$, 
$$ R^{\Quant_k}(\eta, \mu)_ j = i k^{-1} \Quant_k (A (f,g))_j + O(k^{-2}).$$
and $\Quant_k (A (f,g))_j$ is not $O(k^{-1})$ as soon as $A(f,g)$ doesn't vanish.

The first equation of the theorem can
be deduced from the expression of the curvature in theorem
\ref{theo:courbureP} as follows. First recall that 
$$ ik \bigl[ f + \tfrac{1}{ik} \nabla_{X}^{L^k} , g + \tfrac{1}{ik} \nabla_{Y}^{L^k} \bigr] =
\{ f,g \} + \tfrac{1}{ik} \nabla_{[X,Y]}^{L^k},$$
where $X$ and $Y$ are the Hamiltonian vector fields of $f$ and $g$.  
Then we compute the bracket of the operators $D_X$, $D_Y$ entering in
the definition of $\Op_k(f)$ and $\Op_k(g)$ (cf. \eqref{eq:defDX}) in terms of the infinitesimal
actions $\eta$ and $\mu$ on $\Comp$ of $X$ and $Y$ respectively.
\begin{lemme} 
We have 
$ \bigl[ D_X , D_Y \bigr] = D_{[X,Y]} +  \frac{1}{2} \trace(
  \bar{\mu} \eta -  \bar{\eta}\mu )$
\end{lemme} 
\begin{proof} 
Given any one-form $\be$ and complex structure $j$, denote by
$p^{1,0}_j \be$  and $p^{0,1}_j \be$ the component of $\be$ of type
$(1,0)$ and $(0,1)$ for $j$. 
The vector field $\eta$ is given at $j$ by 
$$ \eta_j \in \Ci (M,\Hom( \La_{j}^{1,0}T^*M , \La_{j}^{0,1} T^*M )), \qquad
\eta_j (\al ) = p^{0,1}_j \Lie_X \al.   $$
Consequently, if $\al \in
\Om_j^{1,0} M $ 
\begin{xalignat*}{2} 
[ p_j^{1,0} \Lie_X, p_j^{1,0} \Lie_Y ] \al & = p_j^{1,0} \Lie_X (
\Lie_Y - p^{0,1}_j \Lie_Y ) \al - p_j^{1,0} \Lie_Y (
\Lie_X - p^{0,1}_j \Lie_X ) \al
\\ & = (p_j^{1,0} \Lie_{[X,Y]} - \bar{\eta}_j \mu_j + \bar{\mu}_j \eta_j ) \al
\end{xalignat*} 
So if $p_j$ is the projection from $\La^nT^*M \otimes \C$ onto
$\La_{j}^{n,0}T^*M $ with kernel $\La_{j}^{n-1,1}T^*M \oplus ...\oplus
\La_{j}^{0,n}T^*M $, we have for any $(n,0)$-form $\al$
$$ [ p_j \Lie_X, p_j \Lie_Y ] \al = (p_j \Lie_{[X,Y]} +
\trace (\bar{\mu}_j \eta_j -  \bar{\eta}_j \mu_j) ) \al  $$
which implies the result.
\end{proof} 
Consequently, 
$$  ik [ \Op_k (f) , \Op_k (g) ] = \Op_k ( \{ f,g \} ) +\tfrac{1}{2ik} \trace(
  \bar{\mu}\eta -  \bar{\eta}\mu )$$
and we deduce the first equation of theorem \ref{theo:courbcom} from theorem \ref{theo:courbureP}.

To prove the second equation, we start with the following relation
between the curvatures of
$\Preq_k$ and $\Quant_k$.
\begin{lemme} \label{lem:imp}
For every vector fields   $\eta$ and $\mu$ 
of $\Compint$, we have
$$ R^{\Quant_k}(\eta, \mu) = \Pi_k \bigl[  \nabla_{\eta} ^{
  \End(\Preq_k)} \Pi_k ,  \nabla_{\mu} ^{
  \End(\Preq_k)} \Pi_k  \bigr] + \Pi_k R^{\Preq_k} ( \eta, \mu)$$
where $\nabla^{ \End(\Preq_k)} \Pi_k$ is the commutator $ [ \nabla^{ \Preq_k}, \Pi_k].$
\end{lemme} 
\begin{proof} 
We have 
\begin{xalignat*}{2}
\nabla_{ \eta}^{\Quant_k}  \nabla_{\mu} ^{ \Quant_k} = & \Pi_k
\nabla_{ \eta}^{\Preq_k} \Pi_k \nabla_{\mu} ^{\Preq_k} \Pi_k \\
= & \Pi_k \bigl( \nabla_{ \eta} ^{ \End(\Preq_k) } \Pi_k \bigr)
\nabla_{\mu} ^{\Preq_k} \Pi_k + \Pi_k
\nabla_{ \eta}^{\Preq_k}  \nabla_{\mu} ^{\Preq_k} \Pi_k \\
= &\Pi_k \bigl(  \nabla_{\eta} ^{
  \End(\Preq_k)} \Pi_k\bigr) \bigl(    \nabla_{\mu} ^{
  \End(\Preq_k)} \Pi_k \bigr)   + \Pi_k  \bigl(  \nabla_{\eta} ^{
  \End(\Preq_k)} \Pi_k\bigr) \Pi_k   \nabla_{\mu} ^{\Preq_k}   \\
 &+ \Pi_k
\nabla_{\eta} ^{ \Preq_k} \nabla_{\mu} ^{ \Preq_k} \Pi_k
\end{xalignat*}  
Since $\Pi_k^2 = \Pi_k$, we have $ \Pi_k \bigl( \nabla ^{\End ( \Preq_k) }
\Pi_k \bigr) \Pi_k =0$. So the second term of the sum vanishes. Using this
it is easy to compute the curvature of $\Quant_k$ 
\begin{xalignat*}{2} 
  R^{\Quant_k}(\eta, \mu) & = [ \nabla_{
  \eta}^{\Quant_k},  \nabla_{\mu} ^{ \Quant_k}] - \nabla_{[
  \eta, \mu ]}^{\Quant_k} \\ 
& = \Pi_k \bigl[  \nabla_{\eta} ^{
  \End(\Preq_k)} \Pi_k ,  \nabla_{\mu} ^{
  \End(\Preq_k)} \Pi_k  \bigr] +  \Pi_k
[\nabla_{\eta} ^{ \Preq_k}, \nabla_{\mu} ^{ \Preq_k}]  - \Pi_k \nabla_{[
  \eta, \mu ]}^{\Preq_k} \\ 
& = \Pi_k \bigl[  \nabla_{\eta} ^{
  \End(\Preq_k)} \Pi_k ,  \nabla_{\mu} ^{
  \End(\Preq_k)} \Pi_k  \bigr] +  \Pi_k R^{\Preq_k} ( \eta, \mu)
\end{xalignat*} 
which proves the result. 
\end{proof} 

On the other hand we can compute the commutator of $\Pi_k$ with $\Op_k(f)$
in terms of the covariant derivative of $\Pi_k$. 

\begin{lemme} Let $f \in \Ci (M, \R)$ and $\eta$ be the infinitesimal
  action of $f$ on $\Comp$, then 
$$ \tfrac{1}{ik} \nabla_{\eta}^{ \End(\Preq_k)} \Pi_k = \bigl[ \Op_k (f), \Pi_k \bigr] .$$
\end{lemme} 

\begin{proof} 
This follows from theorem \ref{theo:infinit} by derivating the relation
$\Pi_{k,j_t} U_t = U_t \Pi_{k,j_0}$. 
\end{proof}

Applying twice this last lemma, we obtain
\begin{xalignat*}{2}
\Pi_k \Op_k (f) \Pi_k \Op_k (g) \Pi_k = & \Pi_k \Op_k (f) \Pi_k ( \Op_k (g) + \tfrac{1}{ik}
\nabla_{\mu}^{ \End(\Preq_k)} \Pi_k) \\
 = & \Pi_k ( \Op_k (f)  + \tfrac{1}{ik}
\nabla_{\eta}^{ \End(\Preq_k)} \Pi_k ) ( \Op_k (g) + \tfrac{1}{ik}
\nabla_{\mu}^{ \End(\Preq_k)} \Pi_k) 
\end{xalignat*}
Hence
\begin{gather} \label{eq:eq1} 
\bigl[ \Pi_k \Op_k (f) \Pi_k, \Pi_k \Op_k (g) \Pi_k \bigr] = \\ \notag
\Pi_k \bigl[   \Op_k (f)  + \tfrac{1}{ik}
\nabla_{\eta}^{ \End(\Preq_k)} \Pi_k ,\ \Op_k (g) + \tfrac{1}{ik}
\nabla_{\mu}^{ \End(\Preq_k)} \Pi_k \bigr] \end{gather}
Similarly, we have
\begin{xalignat*}{2}
\Pi_k \Op_k (f) \Pi_k \Op_k (g)\Pi_k  = & ( \Op_k (f) - \tfrac{1}{ik}
\nabla_{\eta}^{ \End(\Preq_k)} \Pi_k ) \Pi_k  \Op_k (g) \Pi_k  \\
 = & ( \Op_k (f)  - \tfrac{1}{ik}
\nabla_{\eta} ^{ \End(\Preq_k)}\Pi_k ) ( \Op_k (g) - \tfrac{1}{ik}
\nabla_{\mu}^{ \End(\Preq_k)} \Pi_k) \Pi_k
\end{xalignat*}
So
\begin{gather} \label{eq:eq2} 
\bigl[ \Pi_k \Op_k (f) \Pi_k, \Pi_k \Op_k (g) \Pi_k \bigr] = \\ \notag \bigl[   \Op_k (f)  - \tfrac{1}{ik}
\nabla_{\eta}^{ \End(\Preq_k)} \Pi_k ,\ \Op_k (g) - \tfrac{1}{ik}
\nabla_{\mu}^{ \End(\Preq_k)} \Pi_k \bigr]\Pi_k \end{gather}
Now equations \eqref{eq:eq1} and \eqref{eq:eq2} imply
\begin{gather*}  
\bigl[ \Pi_k \Op_k (f) \Pi_k, \Pi_k \Op_k (g) \Pi_k \bigr] = \\
 \Pi_k \bigl[ \Op_k
(f), \Op_k (g) \bigr] \Pi_k + \Pi_k \bigl[ \tfrac{1}{ik}
\nabla_{\eta}^{ \End(\Preq_k)} \Pi_k ,\ \tfrac{1}{ik}
\nabla_{\mu} ^{ \End(\Preq_k)}\Pi_k \bigr] \Pi_k 
\end{gather*}
And we deduce the second equation of theorem \ref{theo:courbcom} from the first one and
lemma \ref{lem:imp}.  

\subsection{An analog result in finite dimension} \label{sec:fdr}

It is interesting to note that the expression for the curvature in
theorem \ref{theo:courbureP} is a direct
consequence at least formally of theorem \ref{theo:infinit} on the
infinitesimal action.  Consider a finite dimensional
vector bundle $E \rightarrow X$ endowed with a connection
$\nabla$. Assume a Lie group $G$ acts on $E$ preserving the connection. Given $\eta$ in the Lie
algebra $\mathfrak{g}$ of $G$, denote by $\eta_X$ the vector field
corresponding to the infinitesimal action on the base $X$ and by
$\Lie_{\eta}$ the infinitesimal action on $\Ci( X, E)$. Then
$\Lie_{\eta} - \nabla_{\eta_X}$ acts by exterior multiplication by a
section 
 $$M (\eta) \in \Ci ( M, \End(E)) $$
called the moment of $\eta$. 

\begin{prop} 
For any vectors $\eta, \mu \in {\mathfrak{g}}$, we have
$$ [ M(\eta), M( \mu) ] = M( [\eta, \mu ] ) + R^E (\eta_X, \mu_X),$$
where $R^E$ is curvature of $\nabla$.
\end{prop}

\begin{proof} Since the connection in invariant, we
have
$$ [ \Lie_{\eta}, \nabla ] =0$$
Replacing $\Lie_{\eta}$ with $\nabla_{\eta_X} + M (\eta)$, we obtain that
\begin{gather}  \label{eq:eq}
\nabla^{\End (E)} M (\eta) = R^E (\eta_X, .). \end{gather}
 Since $\eta \rightarrow \nabla_{\eta_{X}} + M(\eta)$ is a Lie algebra
representation, we have for any $\eta, \mu \in
\mathfrak{g}$
\begin{xalignat*}{2}  
 [ \nabla_{\eta_X} + M (\eta) ,\nabla_{\mu_X} + M (\mu) ] = &
\nabla_{[\eta, \mu]_{X} } + M ([\eta, \mu]) \\
= &
\nabla_{[\eta_X, \mu_X]} + M ([\eta, \mu])
\end{xalignat*} 
Assuming that $\nabla$ is $G$-invariant, we can compute by the lemma
the commutators 
$$ [ \nabla_{\eta_X} , M(\mu ) ] = R^E ( \mu_X, \eta_X), \qquad [M (\eta)
,\nabla_{\mu_X}]= - R^E ( \eta_X,  \mu_X) = R^E ( \mu_X, \eta_X).$$
Then using that 
$$  [ \nabla_{\eta_X} ,\nabla_{\mu_X} ] = \nabla_{[\eta_X, \mu_X]}
+R^E ( \eta_X, \mu_X)$$
we obtain the proposition. \end{proof} 

If we apply this equation in our infinite dimensional setting with
$E$ the bundle of quantum spaces or prequantum spaces and $G$ the group of
prequantization bundle automorphism, we
obtain theorem \ref{theo:courbcom}.

It is also interesting to consider the situation the introduction
without half-form bundle (cf. section \ref{sec:intw}). Let $$\pqs_k:= \Ci
(M, L^k) \times \Compint$$ be
the prequantum space bundle and $\qs_k$ be the subbundle of quantum
spaces. As explained in the introduction the Group $\Group$ of
prequantization bundle automorphisms acts on these bundle. Moreover
these bundles are endowed with invariant connection. Then one proves that the moment of
$f \in \Ci (M)$ on $\qs_k$ is the operator $(ik) \qs_k(f)$,
where $\qs_k (f)$ is the Toeplitz operator
$$ \qs_k (f) := \Pi_k \bigl( f + \frac{1}{ik} \nabla^{L^k}_{X} \bigr)
$$
with $X$ the Hamiltonian vector field of $f$. Consequently, one has 
$$ ik [ \qs_k (f) , \qs_k (g) ] = \qs_k ( \{ f,g \} ) + (ik)^{-1} R^{\qs_k}(\eta, \mu)$$
where $\eta$ and $\mu$ are the infinitesimal actions of $f$ and $g$
respectively on $\Compint$. Then we recover the
main point of the argument of Ginzburg and Montgomery: if the
curvature vanishes, the map 
$$ f \rightarrow (ik) \qs_k(f)$$ is a Lie algebra
representation. Furthermore, the result of Foth and Uribe gives the
first correction terms in the computation of the commutator of two
Toeplitz operators.

\section{Preliminaries for the proof of theorem \ref{theo:courbureQ}}

Given two complex structures $j_a, j_b$, we introduce a class of
operators from $\Preq_{k,j_a}$ to $\Preq_{k,j_b}$ extending the class of
Fourier integral operator we considered previously. First by using the
scalar product of $\Preq_{k,j_a}$, the Schwartz kernels of these
operators can be regarded as $\Ci$ sections of the bundle  
$$ \bigl( L^k \otimes \delta_{j_b} \bigr) \boxtimes \bigl( \bar{L}^{k}
\otimes \bar{\delta}_{j_a} \bigr) \rightarrow M^2. $$
Let $N$ be a non-negative integer. We say that $(T_k)_{k \in \N^*}$ is an operator of $\algebre_N ( j_a ,
j_b)$ if its Schwartz kernel is of the form
$$ T_k(x,y) =  \Bigl( \frac{k}{2\pi} \Bigr)^{n} E^k(x,y) f(x,y,k) + O
(k^{-\infty}) $$
where 
\begin{itemize}
\item 
$E$ is a section of  $L \boxtimes \bar{L}
\rightarrow M^2$ such that  $\| E(x,y)
\| <1 $ if $x \neq y$, 
$$ E (x,x) = u \otimes \bar{u}, \quad \forall u \in L_x \text{ such that }
  \| u \| = 1, $$
and $ \bar{\partial} E \equiv 0 $
modulo a section vanishing to any order along $\De$.
\item
  $f(.,k)$ is a sequence of sections of  $ \delta_{j_b}   \boxtimes
  \bar{\delta}_{j_a}
\rightarrow V$ which has an asymptotic expansion in the $\Ci$ topology
  
$$ f(.,k) = k^{N}f_{-N}  + k^{N-1} f_{-N+1} + ...$$
where for all $ 0\leqslant l \leqslant
N$, $f_{-l}$ vanishes to order $2l$ along the diagonal of $M^2$.  
\end{itemize}

As a result, the Schwartz
kernel of $T_k$ is uniformly $O(k^{n +N})$. It is  $O(k^{n +N
  -\frac{1}{2}})$ if and only if $f_l$ vanishes at order $2l +1$ along the
diagonal whenever $2l +1 \geqslant 0$. This follows from the fact that $\ln \| E \| <0$ outside the diagonal and its Hessian along the
diagonal is non-degenerate in the transverse directions
(cf. lemma 1 in \cite{oim1}). We define
the symbol of $(T_k)$ as 
$$ \hb^{-N} [f_{-N} ]_{2N} + \hb^{-N+1} [f_{-N+1} ]_{2(N-1)} +....+\hb^{-1}
[f_{-1}]_2 + [f_0] $$
where $ [f_{-N+l}]_{2(N-l)}$ is the equivalence class of $f_{N-l}$ modulo
the functions vanishing at order $2(N-l) +1$ along the diagonal. So with
the usual identification, the space of symbol is the space of sections
of 
$$  S_N^{j_a,j_b} := \delta_{j_b} \otimes  \bar{\delta}_{j_a} \otimes
\Bigl( \hb^{-N} \Sym_{2N} {\mathcal{C}} \oplus  \hb^{-N+1 } \Sym_{2(N-1)}
{\mathcal{C}}    \oplus... \oplus \Sym_0  {\mathcal{C}}  \Bigr)
$$
where ${\mathcal{C}}$ is the conormal bundle of
the diagonal of $M^2$. 

\begin{theo} \label{theo:comp_alg}
The composition of $S \in \algebre_N (j_b,j_c)$ with  $S' \in
\algebre_{N'} (j_a,j_b)$ is an operator of $\algebre_{N + N'}
(j_a,j_c)$. Furthermore there exists a bilinear bundle map 
$$ L^{j_a,j_b,j_c}_{N, N'} : S^{j_b,j_c}_N \times S^{j_a,j_b}_{N'} \rightarrow S^{j_a,j_c}_{N + N'}$$
such that the principal symbol of $S S'$ is $L_{N, N'}^{j_a,j_b,j_c}
(\sigma,\sigma')$ if $\sigma$ and $\sigma '$ are the symbols of $S$
and $S'$ respectively.
\end{theo} 

\begin{proof} The proof is essentially the same as the one of
theorem \ref{P1} about the composition of Fourier integral
operators. The computation of the symbol follows from the version of
the stationary phase lemma stated in the appendix of \cite{firstpart}. \end{proof} 

By extending the Fourier integral operators of $\Fourier ((j_a,
\delta_{j_a}),(j_b, \delta_{j_b}))$ to
operators $\Preq_{k,j_a} \rightarrow \Preq_{k,j_b}$ in such a way that
they satisfy $$\Pi_{k,j_a} T_k \Pi_{k,j_b} = T_k, \quad k=1,2,...$$ 
$\Fourier ((j_a,\delta_{j_a}),(j_b, \delta_{j_b}))$ becomes a subspace of 
$\algebre_0 (j_a,j_b)$. Both definitions of principal symbols are the
same if we identify the sections of $S_0^{j_a,j_b} =\delta_{j_b} \otimes
  \bar{\delta}_{j_a}$ with the fiber bundle morphisms $\delta_{j_a} \rightarrow
  \delta_{j_b}$ by using the scalar product of $\delta_{j_a}$. 

\begin{theo} \label{theo:alg_FIO}
If $T$ is an operator of $\algebre_N(j_a,j_b) $ with symbol $\si$,
then $\Pi_{k,j_a}  T \Pi_{k,j_b}$ is a Fourier integral operator of $\Fourier
((j_a, \delta_{j_a}),(j_b, \delta_{j_b}))$. Furthermore the symbol of $\Pi_{k,j_a} T \Pi_{k,j_b}$ is $
L_N^{j_a,j_b} (\si)$, where  
$$L_N^{j_a,j_b}: S^{j_a,j_b}_N \rightarrow S^{j_a,j_b}_0$$ is a fiber-bundle morphism.  
\end{theo}

\begin{proof} Again the proof relies on the methods of section \ref{sec:preuve_P1}. To show that $\Pi_{k,j_a} S \Pi_{k,j_b}$ is a Fourier integral
  operator, we argue as in the following of lemma \ref{lem:comp_phas}. The
  other part is an application of  the stationary phase lemma in the
  appendix of \cite{firstpart}.
\end{proof}

Finally let us describe explicitly the symbol product for the
composition of operators of $\algebre_N (j_a,j_a)$. To do this it is
convenient to introduce complex coordinates $(U, z^i)$ for $j_a$ and
write the symbol in the following way
$$ \si ( \hb, \bar{Z}, Z, x) = \sum_{l=0}^{N} \hb^{-l} \si_l( \bar{Z}, Z ) (x),
\qquad  x\in U $$
with 
$$ \si_l (\bar{Z}, Z)(x)  =  \sum_{| \al |+ | \be | = 2l }
  \frac{1}{\al ! \be!} \Bigl( \nabla_{(
  \partial_{\bar{z}^1}, 0)} ^{\al(1)}...\nabla_{(
  \partial_{\bar{z}^n}, 0)} ^{\al(n)} \ \nabla_{(0,
  \partial_{z^1})}^{\be(1)} ... \nabla_{(0,
  \partial_{z^n})}^{\be(n)} f_l\Bigr) (x,x) \ \bar{Z}^{\al} Z^\be 
$$
and  $\nabla$ a covariant derivation of $\delta_{j_a} \boxtimes
\bar{\delta}_{j_a}$. 

\begin{theo} \label{theo:comp_form}
With the previous notations, the map of theorem \ref{theo:comp_alg} is 
$$ L_{N, N'}^{a,a,a} ( \si , \si') ( \hb, \bar{Z}, Z, x) =
\sum_{l=0}^{N+N'} \frac{\hb^{l}}{l!}  \Bigl[ \De^l \bigl(  \si  ( \hb, \bar{Z} -
\bar{Y} , Y, x). \si'( \hb, \bar{Y}, Z - Y , x) \bigr) \Bigr]_{\bar{Y} = Y =
  0} $$
where $\De$ is the operator  
$$ \De := \sum_{i,j}  G^{i,j}(x) \partial_{Y^i} \partial_{\bar{Y}^j} $$ 
with $(G^{i,j})$ the inverse matrix of $(G_{j,i})$ whose coefficients
are such that $\om = i \sum $ $G_{i,j} dz^i \wedge d \bar{z}^j$. 
\end{theo}
\scratch{Le produit de $\si$ et $\si'$ que j'utilise est la
  contraction des deux facteurs du milieu de  $\si(x) \otimes \si' (x)
  \in \delta_x \otimes \bar{\delta}_x \otimes \delta_x \otimes
  \bar{\delta}_x$. On peut aussi identifier $\si$ et $\si'$ {\`a} des
  fonctions et faire un b{\^e}te produit.}

There isn't any difficulty to extend these results to the case where
the complex structures depends smoothly on a parameter. 
We end these preliminaries with the variations of the
section $E$ as a function of the complex structure. Let $x \in
M$ and $\Ga$ be a germ at $x$ of a Lagrangian sub\-ma\-nifold of
$M$. Let us fix a unitary section $s$ of $L\rightarrow \Ga$. Let $j_t$
 be a curve in $\Compint$. Then consider a smooth family $E_t$ of
sections of $L \rightarrow M$ such that $ E_t = s$ along $\Ga$  
and  
$$ \bar{\partial}_{j_t} E_t \equiv 0$$
modulo a section vanishing to any order along 
$\Ga$. Let us write 
$$ \frac{d}{dt} E_t = f_t E_t.$$ 
on a neighborhood of  $x$.
\begin{prop} \label{prop:varE} 
The function $f_0$ and its first derivatives vanish over
$\Ga$. Furthermore, if $Z$ and $Z'$ are holomorphic vector fields
for the complex structure $j_0$, then 
$$ \bar{Z}. \bar{Z}'. f_0 = \frac{1}{i} \om ( \bar{Z}, \mu(\bar{Z}')) 
$$
along $\Ga$, with $\mu \in \Ci (M, \Hom(\La^{1,0}_{j_0}T^*M ,
\La^{0,1}_{j_0} T^*M) )$ the tangent vector to $j_t$ at $j_0$. 
\end{prop}  
Since $\Ga$ is Lagrangian, $T_{j_t}^{0,1} M$ and $T\Ga \otimes \C$ are
transverse. So the result gives the Hessian of $f$ along $\Ga$. 

\begin{proof} 
Since $E_t (y)$ is constant for all $y \in \Ga$, $f_t$
vanishes along $\Ga$. Let us write
$$\nabla E_t = \tfrac{1}{i} E_t \otimes \al_t. $$ 
We can prove that $\al_t$ vanishes along $\Ga$ and  
\begin{gather} \label{eq:noth}
 \bar{Z} . \langle \al_0, W\rangle = \om( \bar{Z},W)\end{gather}
 if $Z$ is a holomorphic vector fields
for $j_0$ (cf. Lemma 4.2 in \cite{firstpart}). 
Since  $\frac{d}{dt}$ et $\nabla$ commute, we have $ d f_t =
\frac{1}{i} \dot{\al}_t$. So $df_t$ vanishes along $\Ga$ because the
 same holds for $\alpha_t$.  

Now let us associate to $j_t \in \Comp$ the section $\mu_t$
 of $\Hom( \La_{j_0}^{1,0}T^*M , \La_{j_0}^{0,1}T^*M)$ as we did in \eqref{eq:defmu}.
Then $T^{0,1}_{j_t}M$ is the graph of 
$$ - \mu_t^t: T_{j_0}^{0,1}M
\rightarrow T_{j_0}^{1,0}M.$$ 
So if $Z'$ is a section of $T_{j_0}^{1,0}M$,  $\bar{Z}'  -
\mu_t^t(\bar{Z}')$ is a section of $  T_{j_t}^{0,1}M$ and  
$$\langle \al_t , \bar{Z}'  - \mu_t^t(\bar{Z}')\rangle $$ 
vanishes to any order along $\Ga$. Thus the same holds for the
derivative 
$$ \langle \dot{\al}_t , \bar{Z}'  - \mu_t^t(\bar{Z}')\rangle -
\langle \al_t ,\dot{\mu}_t^t(\bar{Z}')\rangle .$$
In particular, we have along $\Ga$,  
$$ \bar{Z} .\langle \dot{\al}_0 , \bar{Z}' \rangle = \bar{Z}. \langle \al_0
,\dot{\mu}_0^t(\bar{Z}')\rangle $$
Here we used that $\mu_0 =0$. Finally the results follows from $ d f_t =
\frac{1}{i} \dot{\al}_t$ and \eqref{eq:noth}.
\end{proof}

\section{Proof of theorem \ref{theo:courbureQ}}  \label{sec:preuve}

Let us start with the computation of the curvature.  By lemma
\ref{lem:imp}, 
$$ R^{\Quant_k}(\eta, \mu) = \Pi_k \bigl[  \nabla_{\eta} ^{
  \End(\Preq_k)} \Pi_k ,  \nabla_{\mu} ^{
  \End(\Preq_k)} \Pi_k  \bigr] + \Pi_k R^{\Preq_k} ( \eta, \mu).$$
By theorem \ref{theo:courbureP}, $ \Pi_k R^{\Preq_k} (
\eta, \mu)$ is at $j$ a Toeplitz operator of $\Quant_j$ with
principal symbol 
$\frac{1}{2} \trace(
  \eta.\bar{\mu} - \mu. \bar{\eta} ) (j).$
Then it follows from the following proposition that $R^{\Quant_k}(\eta,
  \mu) (j)$ is a Toeplitz operator with vanishing principal symbol.

\begin{prop} For any tangent vector $\eta, \mu \in T_j \Compint$,
  the operator $$\Pi_k \bigl[  \nabla_{\eta} ^{
  \End(\Preq_k)} \Pi_k ,  \nabla_{\mu} ^{
  \End(\Preq_k)} \Pi_k  \bigr] $$ is a Toeplitz operator of
$\Quant_j$ with principal symbol  $- \frac{1}{2} \trace(
  \eta.\bar{\mu} - \mu. \bar{\eta} )  (j)$.
\end{prop} 

\begin{proof} First we prove that   $\nabla_{\eta} ^{
  \End(\Preq_k)} \Pi_k$ is an operator of 
$\algebre_{2}(j,j)$ and compute its symbol. Let $j_t$ be a curve of
$\Compint$ whose tangent vector at $0$ is $\eta$. Let $\Psi_{t} : \delta_{j_0} \rightarrow
\delta_{j_t}$ be the continuous family of half-form bundle morphisms
such that $\Psi_0$ is the identity of $\delta_{j_0}$. 
We have at $j_0$
$$  \nabla_{\eta} ^{
  \End(\Preq_k)} \Pi_k . \Phi  =   \frac{d}{dt} \Bigr|_{t=0} \Psi^{-1}_{t} \Pi_{k,j_t}
\Psi_{t} \ \Phi   , \qquad \Phi \in \Preq_{k,j_0} $$
\scratch{en effet, comme $\nabla_{\mu} ^{
  \End(\Preq_k)} \Pi_k $ agit par une application lin{\'e}aire dans chaque
fibre, on peut prolonger $\Phi \in \Preq_{k,j_0}$ arbitrairement et
calculer $$\nabla_{\mu} ^{
  \End(\Preq_k)} \Pi_k . \Phi = \bigr( \nabla_{\mu} ^{\Preq_k} \circ \Pi_k - \Pi_k \circ
\nabla_{\mu} ^{\Preq_k} \bigl) . \Phi$$
en $j_0$. On s'arrange pour que $\Phi (j_t)=\Psi_t . \Phi(j_0)$ et  {\c c}a donne la formule pr{\'e}c{\'e}dente (cette section sera
plate en $j_0$).}
Recall that $ \Pi_{k,j_t}$ is an operator of  $\algebre_0
(j_t,j_t)$ with symbol 1. Thus its kernel is of the form  
$$  \Bigl( \frac{k}{2\pi} \Bigr)^{n} E_t^k(x,y) f_t(x,y,k) + O
(k^{-\infty}) $$
where  $f_t(.,k)$ is a sequence of sections of $\delta_{j_t}
\boxtimes \bar{\delta}_{j_t}$ equal to $1+O(k^{-1})$ on the
diagonal. We obtain the kernel of $\Psi^{-1}_{t} \Pi_{k,j_t}
\Psi_{t}$ by replacing  $f_t$ with  $$\Psi^{-1}_{t} (x) f_t(x,y,k)
\Psi_{t} (y)$$ which again is equal to  $1  +O(k^{-1})$ on the
diagonal. Derivating with respect to $t$, we deduce that $\nabla_{\eta} ^{
  \End(\Preq_k)} \Pi_k$ is an operator of 
$\algebre_{2}(j,j)$ with symbol
$$ \hbar^{-1} [g]$$
where  $g$ is the function of $M^2$ such that  $\frac{d}{dt} \Bigr|_{t=0} E_t = g
E_0$. 

Let us compute $[g]$. Let $(z^i)$
be a complex coordinates system for $j_0$ such that $\om =
i \sum dz^j \wedge d\bar{z}^j$ at $x$. Denote by $U(t)$ 
the symmetric matrix such that the family  
$$ dz^i + \textstyle{\sum}_j U_{ij}(t) d\bar{z}^j, \qquad i =1,...,n$$
is a base of $\La_{j_t}^{1,0} T_x^*M$. So the derivative $\dot{U}$
of $U(t)$ at $t=0$ is the matrix of $\eta (j_0)$. By proposition
\ref{prop:varE}, the symbol of $\nabla_{\eta} ^{
  \End(\Preq_k)} \Pi_k$ at $x$ is   
$$ \hbar^{-1} [g](\bar{Z}, Z,x) = - \frac{1}{2\hb} \sum ( \dot{U}_{ij} \bar{Z}^i \bar{Z}^j +
\dot{\bar{U}}_{ij} Z^i Z^j) $$
where we used the notations of  theorem \ref{theo:comp_form}. 

Then it follows from theorems \ref{theo:comp_alg} and \ref{theo:alg_FIO} that $$\Pi_k \bigl[  \nabla_{\eta} ^{
  \End(\Preq_k)} \Pi_k ,  \nabla_{\mu} ^{
  \End(\Preq_k)} \Pi_k  \bigr]$$ is a Toeplitz operator and we compute its
symbol by applying theorem \ref{theo:comp_form}. 
First the symbol  of $\bigl[  \nabla_{\eta} ^{
  \End(\Preq_k)} \Pi_k ,  \nabla_{\mu} ^{
  \End(\Preq_k)} \Pi_k  \bigr] $ is at $x$
$$ \hb^{-1} \sum \bigl( \dot{U}_{ik} \dot{\bar{V}}_{jk} -\dot{V}_{ik}
\dot{\bar{U}}_{jk} \bigr)   \bar{Z}^i Z^j - \frac{\hb^{-2}}{4} \sum  \bigl( \dot{U}_{ij} \dot{\bar{V}}_{kl} -\dot{V}_{ij}
\dot{\bar{U}}_{kl} \bigr)   \bar{Z}^i \bar{Z}^j Z^k Z^l $$
where $\dot{V}$ is associated to $\mu$ as $\dot{U}$
to $\eta$. Then the symbol of $\Pi_k \bigl[  \nabla_{\eta} ^{
  \End(\Preq_k)} \Pi_k ,  \nabla_{\mu} ^{
  \End(\Preq_k)} \Pi_k  \bigr]$ is at $x$
$$- \frac{1}{2} \trace\ ( \dot{U} \dot{\bar{V}} - \dot{V} \dot{\bar{U}}  ) $$
Since $\dot{U}$ and $\dot{V}$ are the matrices of $\eta(j_0)$ and
$\mu (j_0)$, we obtain the result.
\end{proof}

Let us prove now the last part of theorem
\ref{theo:courbureQ}. Consider a curve $j:[0,1] \rightarrow \Compint$. Denote by $\Psi_t
: \delta_{j_0}
\rightarrow \delta_{j_t}$
the continuous family of half-form bundle morphisms such that $\Psi_0$
is the identity. 

\begin{prop} Consider a smooth family of Fourier integral operator
  $$\bigl( P_t  
\in \Fourier ((j_0,\delta_{j_0}), (j_t, \delta_{j_t})); \;  t\in [0,1]
\Bigr)$$ with symbol $(\sigma_t
\Psi_{t})_t$. Then  the operator
$$ \bigl( \nabla^{j^*\Quant_k}_{\partial_t} \circ P \bigr)_t : \Quant_{k,j_0}
\rightarrow \Quant_{k,j_t}, \quad \Phi_0 \rightarrow
(\nabla^{j^*\Quant_k}_{\partial_t} \Phi)(t) \text{ with } \Phi (t) = P_t
  \Phi_0   $$
is a Fourier integral operator of $\Fourier ((j_0, \delta_{j_0}),(j_t,
\delta_{j_t}))$ and its symbol
is $\dot{\sigma}_t
\Psi_{t}$
\end{prop} 

\begin{proof} The Schwartz kernel of $P_t$ is of the form
$$ \Bigl( \frac{k}{2 \pi} \Bigr)^n E_t^{k} (x,y) f_t (x,y,k) $$
with $f_t (x,x,k) =\si_t (x) \Psi_{t}(x) + O(k^{-1})$. 
So the kernel of $( \nabla^{j^*\Preq_k}_{\partial_t} \circ P )_t$
 is 
$$ \Bigl( \frac{k}{2 \pi} \Bigr)^n  \Bigr[ \frac{d}{dt}  E_t^{k} (x,y)
\Bigl] f_t(x,y,k)+ \Bigl( \frac{k}{2 \pi} \Bigr)^n E^{k}_t (x,y)
\Bigl[ \frac{d}{ds} (\Psi_{t,s} ^{-1} (x). f_s
(x,y,k) )\Bigr]_{s=t} $$
where $(\Psi_{t,s} : \delta_{j_t} \rightarrow \delta_{j_s})_s$ is the
continuous family of half-form bundle morphisms such that $\Psi_{t,t}$
is the identity of $\delta_{j_t}$. By  proposition \ref{prop:varE}, 
$$ \frac{d}{dt}  E_t = g_t E_t, $$ 
where $g_t$ and its first derivatives vanish along the diagonal. So
$( \nabla^{j^*\Preq_k}_{\partial_t} \circ P \bigr)_t$ belongs to $\algebre_2 (j_0,
j_t)$. Since 
$$ \Psi_{j_t, j_s} \circ \Psi_{j_0, j_t} = \bic (j_0, j_t,j_s)
\Psi_{j_0,j_s} $$
we have 
$$ \Psi^{-1}_{t,s} \circ \Psi_{0,s} = \frac{1} {\bi (j_0,j_t,j_s)} \Psi_{0,t}$$
Thus the symbol of $ ( \nabla^{j^*\Preq_k}_{\partial_t} \circ P
\bigr)_t$  is
\begin{gather} \label{eq:symbole}
  \Bigl[  \bigl( \hb^{-1} [g_t] + \frac{d}{ds}\Bigl( \frac{1} {\bi
  (j_0,j_t,j_s)} \Bigr)\Bigr|_{s=t} \bigr) \si_t + \dot{\si}_t \Bigr] \Psi_{0,t}. 
\end{gather}
Then it follows from theorem \ref{theo:alg_FIO} that 
$$ ( \nabla^{j^*\Quant_k}_{\partial_t} \circ P \bigr)_t = \Pi_{k,j_t}
  \circ ( \nabla^{j^*\Quant_k}_{\partial_t} \circ P
  \bigr)_t   $$ belongs to $ \Fourier ((j_0,\delta_{j_0}),(j_t, \delta_{j_t}))$. Furthermore its
  symbol is of the form $ (a_t \si_t+ \dot{\si}_t) \Psi_t$ with $a_t$ a
  $\Ci$ function. To end the proof, it suffices
to show that if $\si_t =1$ for every $t$, then the symbol of $(
  \nabla^{j^*\Quant_k}_{\partial_t} \circ P \bigr)_t$  vanishes. 

Let us check it at $t=0$.  Since $\bi (j_0,j_0,j_s ) =1$ for every $s$
  and
  $\Psi_{0,0}$ is the identity, formula
  \eqref{eq:symbole} simplifies into 
$$   \hb^{-1} [g_0] \si_0 + \dot{\si}_0 $$
which is equal to $\hb^{-1} [g_0] $ because $\si_t =1$. 
Introduce complex coordinates $(z^i)$ for the complex structure $j_0$
such that $\om = i
dz^i \wedge d\bar{z}^i$ at $x$. Let $\mu$ be the tangent vector of
$j_t$ at $t=0$. Let $ U$ be the matrix at $x$ of
$\mu$ in the bases $dz^i$, $d\bar{z}^i$. Then by proposition \ref{prop:varE}
$$ [g_0] = -\frac{1}{2} \sum U_{ij} \bar{Z}^i \bar{Z}^j $$
Finally an application of theorem \ref{theo:comp_form} proves that the symbol of $
\nabla^{\Quant_k}_{\mu_t} \circ P_t$ vanishes at  $t=0$. 

Let us compute now the symbol at any $t$. Since $\si_t =1$, the
operator $P_t$ is invertible with an inverse in $\Fourier ((j_t, \delta_{j_t}),(j_0,\delta_{j_0}))$. So the operator 
$$ P_{t,s} = P_s \circ P^{-1}_t : \Quant_{k,j_t} \rightarrow
\Quant_{k,j_s}$$
belongs to $\Fourier ((j_t, \delta_{j_t}),(j_s,\delta_{j_s}))$ and by theorem \ref{theo:repr}, its
symbol is $\Psi_{t,s}$. It follows from the previous computation that
$$ \bigl( \nabla^{j^*\Quant_k}_{\partial_s} \circ P_{t,.} \bigr)_s$$
belongs to $ \Fourier ((j_t, \delta_{j_t}),(j_s,\delta_{j_s}))$ and its symbol vanishes at
$t=s$. Consequently, the symbol of  
$$ \bigl( \nabla^{j^*\Quant_k}_{\partial_t} \circ P \bigr)_{t} = \bigl(
\nabla^{j^*\Quant_k}_{\partial_s} \circ P_{t,.} \bigr)_{s=t}  \circ P_t $$
vanishes. 
\end{proof}

Then it is easy to construct by successive approximations a smooth family of operators
$$P_t : \Quant(j_0) \rightarrow \Quant (j_t)$$ in 
$\Fourier((j_0, \delta_{j_0}),(j_t,\delta_{j_t}))$ such that $P_0$ is the identity of
$\Quant_{k,j_0}$ and that the total symbol of 
$(\nabla^{j^*\Quant_k}_{\partial_t} \circ P ) _t$ vanishes. Consequently,
$$ \bigl( \nabla^{j^* \Quant_k}_{\partial_t} \circ P \bigr) _t = O(k^{-\infty})$$
where the big $O$ is for the uniform norm of operators and is uniform
with respect to $t$. If 
$T_t: \Quant_{k,j_0} \rightarrow \Quant_{k,j_t}$
is the parallel transport along $j_t$, then 
$$ T_t = P_t - T_t \int_0^t T_{-s} \bigl( \nabla^{j^*
  \Quant_k}_{\partial_s} \circ P \bigr) _s  \ {ds}.$$ 
By the first part of theorem \ref{theo:courbureQ}, $T_t$ is
unitary. Consequently  
$$  T_t = P_t + O(k^{-\infty}).$$
Then using that 
$$\Pi_{k,j_t} (T_t - P_t) \Pi_{k,j_0} = T_t - P_t$$ and $\Pi_{k,j}
\in \Fourier ((j, \delta_{j}),(j,\delta_{j}))$, we show that the Schwartz kernel of 
$T_t -P_t$ is uniformly $O(k^{-\infty})$ with its successive covariant
derivatives. This proves theorem \ref{theo:courbureQ}.

\bibliography{biblio}
\end{document}